%% file: main.tex
\documentclass[a4paper,reqno]{amsart}

\usepackage{xcolor}
\usepackage[utf8]{inputenc}
\usepackage[american]{babel}
\usepackage{csquotes}
\usepackage[T1]{fontenc}
\usepackage{microtype} 
\usepackage{fullpage}
\usepackage{amsmath,amssymb,amsthm}
\usepackage{thmtools}
\usepackage{mathtools} 
\usepackage{upgreek}
\usepackage{orcidlink}
\usepackage{hyperref}
\usepackage{enumitem}
\usepackage{todonotes}
\usepackage[citestyle=numeric, backend=biber, url = false, doi=true, minbibnames=1, maxbibnames=8, defernumbers=false, uniquename=false, giveninits=true]{biblatex}
\allowdisplaybreaks
\hypersetup{hidelinks,hypertexnames=false}
\numberwithin{equation}{section}

\declaretheorem[name=Theorem, numberwithin=section]{theorem}
\declaretheorem[name=Proposition, sibling=theorem]{proposition}
\declaretheorem[name=Lemma, sibling=theorem]{lemma}
\declaretheorem[name=Corollary, sibling=theorem]{corollary}
\declaretheorem[name=Definition, style=definition, sibling=theorem]{definition}

\declaretheorem[name=Example, style=definition, sibling=theorem]{example}
\declaretheorem[name=Remark, style=definition, sibling=theorem]{remark}
\declaretheorem[name=Notation, style=definition, sibling=theorem]{notation}
\declaretheorem[name=Convention, style=definition, sibling=theorem]{convention}

\newcommand{\N}{\mathbb N}
\newcommand{\Z}{\mathbb Z}
\newcommand{\R}{\mathbb R}
\newcommand{\C}{\mathbb C}

\newcommand{\cH}{\mathcal H}

\newcommand{\cD}{\mathcal D}

\newcommand{\cC}{\mathcal C}
\newcommand{\cB}{\mathcal B}

\DeclareMathOperator{\dist}{dist}
\DeclareMathOperator{\rank}{rank}

\DeclareMathOperator{\interior}{int}

\begin{document}

\title{Localizing spectra and pseudospectra\\of a class of unbounded operators on Hilbert space}

\author[A.~Buchinger]{Andreas Buchinger\,\orcidlink{0009-0004-4203-5874}}
\address[A.B.]{TU Hamburg, Institute of Mathematics,
  Am Schwarzenberg-Campus 3, 21073 Hamburg, Germany}
\email{andreas.buchinger@tuhh.de}

\author[M.~Lindner]{Marko Lindner\,\orcidlink{0000-0001-8483-2944}}
\address[M.L.]{TU Hamburg, Institute of Mathematics,
  Am Schwarzenberg-Campus 3, 21073 Hamburg, Germany}
\email{lindner@tuhh.de}

\author[M.~Wittig]{Mattes Wittig\,\orcidlink{0009-0000-3800-390X}}
\address[M.W.]{TU Hamburg, Institute of Mathematics,
  Am Schwarzenberg-Campus 3, 21073 Hamburg, Germany}
\email{mattes.wittig@tuhh.de}

\makeatletter
\hypersetup{
  pdftitle={\@title},
  pdfauthor={A. Buchinger, M. Lindner, M. Wittig},
  pdfkeywords={Keywords}
}
\makeatother

\date{\today}


\begin{abstract}
We use a localization technique to bound and approximate the spectrum and the pseudospectra of a certain class of possibly non-normal, generally unbounded but closed operators on $L^2(X)$, where $X$ is either a) a locally compact second-countable Hausdorff group $G$, not necessarily Abelian, with identity, a proper metric, and a left Haar measure or it is b)~a uniformly discrete subset of $G$ with the counting measure. We treat the discrete and continuous cases b) and a) simultaneously. Our results yield computable subsets and supersets of the pseudospectra of $H$ that converge to the corresponding pseudospectrum in Attouch-Wets distance as the size parameter $L$ in our localization procedure goes to infinity. We illustrate our findings with unbounded band operators on $\ell^2(\Z)$ and with Schr\"odinger operators with $L^\infty_{\rm loc}$ potentials acting on $L^2(\R^d)$.
\end{abstract}

\maketitle

\input{intro}

\input{Pseudospectra_Lower_Norm}

\input{Theorems}

\input{discrete_example} 
\input{schroedinger_example}

\printbibliography

\end{document}

%% file: intro.tex
 \section{Introduction}

The role of spectra in science and engineering can not be overestimated; but in dimensions greater than four, their exact computation is generally impossible, by Abel-Ruffini, whence one is forced to resort to computing subsets, supersets or approximations of the spectrum. So let $H$ be a closed operator on a Hilbert space $\cH$. A typical strategy with a successful history especially for translation invariant operators \cite{BoettcherSilbermann1999,HagenRochSilbermann2001} is to introduce a finite-dimensional projection $P$ on $\cH$ and to study spectral properties of $PHP$ instead of $H$. However, in general, the spectrum of $PHP$ has little to do, even in the limit as $\rank(P)\to\infty$, with the spectrum of $H$.

Our approach, following \cite{ChonchaiyaPhD,ChandlerWildeChonchaiyaLindner2024,HegeMoscolariTeufel2026,LindnerSeifertChandlerWildePreprint,Wittig2026}, is related but a bit different: Roughly speaking, we study all operators of the form $HP$ for a sufficiently large family of projections $P$ (all of the same finite rank but onto different subspaces of $\cH$) and we combine the acquired spectral information in three different ways, where the first provides a subset and the second a superset of the $\varepsilon$-pseudospectrum, $\sigma_\varepsilon(H)$, of $H$, while the third construction yields a superset of the spectrum, $\sigma(H)$, of $H$. We show that all three sets approximate $\sigma_\varepsilon(H)$, resp.~$\sigma(H)$, in an appropriate sense as $\mathrm{rank}(P)\to\infty$.

We want to underline that our projection method-based procedure is {\bf not} ``yet another numerical simulation of what the spectrum $\sigma(H)$ might roughly look like'' -- but it is a method to acquire rigorous subsets and supersets of $\sigma_\varepsilon(H)$ and $\sigma(H)$ -- with error control. Looking, for example, at our superset $S$ of $\sigma(H)$, we can provably exclude points $\lambda\not\in S$ from the spectrum, thereby showing invertibility of $H-\lambda$, and so we can rigorously detect and quantify spectral gaps in infinite-dimensional physical models.

To be a bit more precise, we show that, under certain assumptions on $H$ and the family of all $P$, for every $\varepsilon$-pseudomode\footnote{For $\varepsilon>0$, $\psi\ne0$ is called an $\varepsilon$-pseudomode (or $\varepsilon$-pseudoeigenvector) of $H$ for the pseudoeigenvalue $\lambda\in\C$ if $\|(H-\lambda)\psi\|\le\varepsilon\|\psi\|$. For $\varepsilon=0$ this reduces to the usual meanings of eigenvector and eigenvalue.} $\psi$ of $H$ for the pseudoeigenvalue $\lambda\in\C$, there exists one $P$ in the family such that $P\psi$ is an $(\varepsilon+\delta)$-pseudomode of $HP$ for $\lambda$. 
In accordance with \cite{ChandlerWildeChonchaiyaLindner2024} and \cite{LindnerSeifertChandlerWildePreprint}, let us call this term $\delta$ the truncation penalty. Going through all $P$ in the family, we can hence cover the $\varepsilon$-pseudospectrum of $H$ by the $(\varepsilon+\delta)$-pseudospectra of all $HP$. In contrast, and very simple to argue, the union of the $\varepsilon$-pseudospectra of all $HP$ yields the announced subset of the $\varepsilon$-pseudospectrum of $H$.

The key argument for the superset construction is by averaging both 
\[
\|P(H-\lambda)\psi\|^2 \qquad\text{and}\qquad \|(H-\lambda)P\psi\|^2
\]
over the $P$ family. While the first average is easy to compute, the second requires a detailed analysis of the commutator norm 
\[
\|[H-\lambda,P]\psi\|^2
\]
and its $P$-average. This is essentially the strategy of our forerunners \cite{ChonchaiyaPhD,ChandlerWildeChonchaiyaLindner2024,HegeMoscolariTeufel2026,LindnerSeifertChandlerWildePreprint,Wittig2026} that are all focussed on bounded operators $H$, bounding the average of $\|[H-\lambda,P]\psi\|^2$ by a multiple of $\|\psi\|^2$ -- and the coefficient in this bound is exactly where the truncation penalty $\delta$ comes from.

\begin{samepage}
This explicit expression for the truncation penalty $\delta$ in terms of the commutator average leads to the following conclusions:
\begin{itemize}
\item[(i)] One can show that $\delta\to 0$, and consequently our set sandwich converges into the $\varepsilon$-pseudospectrum of $H$, as the rank of the projections $P$ (being the same for all of them) goes to infinity.
\item[(ii)] One can finetune $\delta$ by passing from projections (sharp cut-offs in $\ell^2$ space) to softer cut-offs $P$, thereby minimizing boundary effects in the commutator and hence the truncation penalty $\delta$.
\end{itemize}
\end{samepage}

This is essentially the common plot of all of \cite{ChonchaiyaPhD,ChandlerWildeChonchaiyaLindner2024,HegeMoscolariTeufel2026,LindnerSeifertChandlerWildePreprint,Wittig2026} -- with some little differences between these papers, of course. So let us first unveil this history and the differences in a bit more detail and then show what our current paper contributes to the story next.

\medskip\noindent {\bf Related work.}
In \cite{ChonchaiyaPhD,ChandlerWildeChonchaiyaLindner2024}, the case of tridiagonal bounded operators $H$ on $\ell^2(\Z)$ was studied with extensions to arbitrary finite bandwidth (via block tridiagonal) and even to the closure of that set. Besides the one-sided truncations, $HP$, also two-sided truncations, $PHP$, and periodized versions of the latter have been studied. \cite{LindnerSeifertChandlerWildePreprint} took two of the three methods to $\ell^2(G)$ for discrete Abelian groups $G$ and, independently and in parallel, \cite{HegeMoscolariTeufel2026} worked out very similar ideas and results for band operators on $\ell^2(\Gamma)$, where $\Gamma$ was a uniformly discrete subset of $\R^d$ that did not have to have any group structure. Then \cite{Wittig2026} came and combined the best of both threads, \cite{LindnerSeifertChandlerWildePreprint} and \cite{HegeMoscolariTeufel2026}, covering the whole story above for band operators on $\ell^2(\Gamma)$, where $\Gamma$ was again a uniformly discrete subset but now of a (not necessarily Abelian) group $G$ with a metric and a Haar measure. One important step in both \cite{HegeMoscolariTeufel2026, Wittig2026} was that, although the vector space $\ell^2(\Gamma)$ was carried by the discrete subset $\Gamma\subset G$, all the localization and averaging procedures were performed over centers in the whole ambient group ($\R^d$ or $G$), so that, in particular, the sums above turn into integrals with the corresponding measure of $\R^d$ or $G$ there.

Besides these immediate forerunners of the current paper, there were several more papers touching our topic: 
For aperiodic discrete Schrödinger operators on $\ell^2(\Z)$, localization techniques for the lower norm even go back to \cite{Elliot1990,AvronMoucheSimon1990}. 
In \cite{LiSei2014,HagLiSei2016}, lower norm localization was shown via the decay of the commutator norm $\|[H-\lambda,P]\|$ as the size of the local neighborhoods in $\ell^2(\Z^d)$ goes to infinity. \cite{SpakulaWillett2016,HaggerSeifert2020} generalize and extend \cite{HagLiSei2016} from $\Z^d$ to metric groups and measure spaces, where \cite{SpakulaWillett2016} is the first explicit account that we know of of an averaged commutator norm analysis as in \cite{LindnerSeifertChandlerWildePreprint,Wittig2026}.
In \cite{LiSchmidt:Givens} it was studied how to numerically benefit from the overlap between the matrices behind $HP$ for different $P$ with overlapping ranges. The approximation of spectra and pseudospectra of unbounded operators by related methods, including the extension of Hausdorff set convergence to this unbounded setting, can be found in many other papers and books like \cite{Hansen2011,B_gli_2014,Boegli_2017, Boegli_2018,Colbrook2026}, to name just a few. Related in a wider sense are of course the many spectral inclusion results by Gershgorin arguments and numerical ranges (of $H$, its polynomials or of its resolvent). We refer to the introduction of \cite{ChandlerWildeChonchaiyaLindner2024} for an extensive review also of that work.

\medskip\noindent {\bf Our contribution.} 
In the current paper we extend that whole plot in two further directions:
\begin{itemize}
\item {\bf to unbounded operators}: The only part of the argument that really requires bounded operators is the analysis of the commutators $[H-\lambda,P]$. If we think of $H$ as represented in $\ell^2(\Z)$ by a bi-infinite matrix $\bf H$ and of $P$ as multiplication by the characteristic function of a $\Z$-interval of length $n$ (with each one of these intervals referring to another projector $P$) then it is not difficult to see that the condition
\[
\sup_P \|[H-\lambda,P]\| < \infty
\]
is equivalent to boundedness of $\bf H-D$, where $\bf D$ refers to the main diagonal of $\bf H$. To give just one example, this relaxation already allows $H$ to be a discrete Schrödinger operator with unbounded potential. But we go much further than this: Instead of bounding the $P$-average of the commutator $\|[H-\lambda,P]\psi\|^2$ by a multiple of the squared $2$-norm, $\|\psi\|^2$, we bound it by a multiple of the squared graph norm, that is,
\begin{equation} \label{eq:comm-graph}
\sum_P \|[H-\lambda,P]\psi\|^2 \le C (\|H\psi\|+\|\psi\|)^2,
\end{equation}
opening the door for a much larger family of unbounded operators $H$.

We argue that the graph norm of $\psi$ is the natural choice of norm at the right-hand side of \eqref{eq:comm-graph} and that it could have appeared at the same position already in the previous papers \cite{ChonchaiyaPhD,ChandlerWildeChonchaiyaLindner2024,HegeMoscolariTeufel2026,LindnerSeifertChandlerWildePreprint,Wittig2026} about bounded operators, where it is however equivalent to the $2$-norm.

\item {\bf from $\ell^2$ to $L^2$}: 
We mentioned already that \cite{HegeMoscolariTeufel2026,Wittig2026} use the ambient group, $\R^d$ or $G$, not the discrete subset $\Gamma$, for the centers of the local neigborhoods and for averaging the norms of $[H-\lambda,P]\psi$, etc. This was already a big step in the direction towards allowing functions on the whole group $G$, thereby moving from operators $H$ on $\ell^2(\Gamma)$ to $L^2(G)$. Because our spectral studies all rest on the evaluation of norms and lower norms, there is indeed not much else to do in order to treat both cases, discrete and continuous, simultaneously in Sections 3 and 4, where $H$ is now an operator on $L^2(X)$ with either a)~$X=G$ equipped with its left Haar measure or b)~$X=\Gamma$ with the counting measure.
\end{itemize}

Because we now have to deal with unbounded spectra and pseudospectra, where Hausdorff convergence is not the right concept anymore, a further contribution of our paper, motivated by \cite{B_gli_2014}, is a study of the Attouch-Wets convergence \cite{Mosco1969,AttouchWets1983,AttouchWets1986} of sequences of closed but not necessarily bounded sets.

\medskip\noindent{\bf Localization: a cost-benefit analysis.}
As always with localization, we translate one tough problem into many simpler problems. So the obvous questions are: a) How much simpler and b) how many?

a) In contrast to the original problem (that often cannot even be stored on a computer), the simpler problems are finite-dimensional\footnote{In the $L^2$-case, finite dimensions are still one spatial discretization away but the reduction to bounded geometry is the key step in that process.} and hence within the reach of numerical linear algebra. More precisely, each of these simpler, local problems is about estimating the smallest singular value of an operator $HP$ on the finite-dimensional space $P(\cH)$. Note that the latter is possible without even actually computing\footnote{See e.g.~Section 7 and in particular (7.4) in \cite{ChandlerWildeChonchaiyaLindner2024}.} that smallest singular value, $\nu(HP)$.

b) The price of localization, besides the truncation penalty $\delta$, is that this family of simpler problems can be pretty large -- and this is when the local structure (or complexity) of $H$ comes into play: Even though the $P$-family is generally infinite, the family of all $HP$ may well be finite and manageable. Citing \cite{HegeMoscolariTeufel2026}, the latter is the case in most physical situations. If we again think of $H$ as represented in $\ell^2(\Z)$ by a bi-infinite matrix $\bf H$ and of each $P$ as multiplication by the characteristic function of another $\Z$-interval of length $n$ then $HP$ is determined by $n$ consecutive columns of $\bf H$. If, as e.g.~in the discrete Schrödinger case, $\bf H$ has only one varying diagonal $\bf v$ then the question is about the set of all subwords of length $n$ in the infinite word $\bf v$. For $p$-periodic $\bf v$, that set has $p$ elements and for aperiodic (1D Sturmian) $\bf v$, it has $n+1$ elements, all of them well-known.

\medskip\noindent {\bf Structure of this paper.}
After this introduction, Section 2 explains open and closed pseudospectra $\sigma_\varepsilon(H)$ and $\Sigma_\varepsilon(H)$ and their study via lower norms before it derives results about the continuity of $\sigma_\varepsilon(H)$ and $\Sigma_\varepsilon(H)$ in terms of both $\varepsilon$ and $H$. A major tool here is a thorough analysis of the Attouch-Wets convergence.
Section 3 explains the geometric setting: the group $G$, the set $X$ and its measure, covering the discrete and the continuous case simultaneously, and the localization functions (the smooth cut-off shapes as mentioned in point (ii) above).
Section 4 is where we actually localize the operators: passing from $H$ to $HP$, study the averages of $\|P(H-\lambda)\psi\|^2$ and $\|(H-\lambda)P\psi\|^2$, in particular that of the commutator, $\|[H-\lambda,P]\psi\|^2$, and draw conclusions about the lower norms, resolvent norms and pseudospectra on a local versus global operator level including the announced set inclusions and convergence results.
Sections~5 and 6 finally illustrate our results with unbounded band operators on $\ell^2(\Z)$ and with Schrödinger operators on $L^2(\R^d)$. In particular, we prove the corresponding commutator average bounds.

%% file: Pseudospectra_Lower_Norm.tex
\section{Spectra, pseudospectra, and lower norms}
\label{sec:spectral-notions}

For a general and in-depth introduction in the theory of pseudospectra, we refer to, e.g., \cite{TrefethenEmbree2005,Colbrook2026}.

We begin by fixing notations used throughout the paper for closed, densely defined operators on a Hilbert-space $\cH$. After recalling the definitions via resolvents and their basic properties, we introduce lower norms as a tool for expressing these notions in a form that is compatible with localisation procedure developed later.  
Let $\cH_1$ and $\cH_2$ be Hilbert spaces. We denote by $\cC(\cH_1,\cH_2)$
the set of closed, densely defined linear operators \[
    T\colon\cD(T)\subset \cH_1 \longrightarrow \cH_2.
    \]
As usual,
$\cB(\cH_1,\cH_2)$ denotes the space of bounded linear operators from
$\cH_1$ to $\cH_2$. If $\cH_1 = \cH_2 = \cH$, we simply write $\cC(\cH)$ and $\cB(\cH)$. Throughout the paper, let $\cH_1,\cH_2,\cH$ always denote Hilbert spaces.
For $H \in \cC(\cH)$,
its resolvent set consists of all
$\lambda\in\mathbb C$ for which\footnote{As is common, we identify $\lambda$ and $\lambda\operatorname{id}_{\cH}$.}
$H-\lambda\colon\mathcal D(H)\longrightarrow\mathcal H$
is bijective with bounded inverse, and its spectrum is the complement of the
resolvent set:
\[
  \rho(H)
  \coloneqq
  \left\{
    \lambda\in\mathbb C:
    (H-\lambda)^{-1}\in\mathcal B(\mathcal H)
  \right\},
  \qquad
  \sigma(H)\coloneqq\mathbb C\setminus\rho(H).
\]
Here and in the following, the bijectivity of $H-\lambda$ is implicitly assumed in the notation
$(H-\lambda)^{-1}$. Since $H$ is closed, boundedness of the inverse actually already follows
from the closed graph theorem once bijectivity is known. Throughout, we use the convention $\frac{1}{0}\coloneqq\infty$ and $\frac{1}{\infty}\coloneqq0$.

\begin{definition}[Pseudospectrum]
\label{def:pseudospectrum}
Let $H\in\cC(\cH)$ and $\varepsilon \geq 0$. We define the \emph{open and closed $\varepsilon$-pseudospectra} of $H$, respectively, by\footnote{For $\varepsilon>0$, one may equivalently omit the explicit union with $\sigma(H)$ by setting
$\|(H-\lambda)^{-1}\|\coloneqq  \infty$ for $\lambda\in\sigma(H)$.}
\begin{align*}
    \sigma_\varepsilon(H)&\coloneqq\sigma(H)\cup\left\{\lambda\in\rho(H):\|(H-\lambda)^{-1}\|>\frac{1}{\varepsilon}\right\}, \\
  \Sigma_\varepsilon(H)&\coloneqq\sigma(H)\cup\left\{\lambda\in\rho(H):\|(H-\lambda)^{-1}\|\geq\frac{1}{\varepsilon}\right\}.  
\end{align*}
For $\varepsilon>0$, these sets are open and closed, respectively. At $\varepsilon=0$, both definitions give
$\sigma_0(H)=\Sigma_0(H)=\sigma(H)$.
\end{definition}

The families $(\sigma_\varepsilon(H))_{\varepsilon\geq0}$ and $(\Sigma_\varepsilon(H))_{\varepsilon\geq0}$ are increasing 
in $\varepsilon$, and and the individual pseudospectra are generally not bounded as $H$ is unbounded. However, if $H\in\mathcal{B}(\mathcal{H})$, then $\sigma_\varepsilon(H)\subset\Sigma_\varepsilon(H)\subset \overline{B_{\lVert H\rVert+\varepsilon}(0)}$.

Although this resolvent-oriented definition is conceptually natural, it is not
directly suited to the localisation arguments used below. In particular, the
local finite sections arising in our setting are generally rectangular and
therefore do not properly correspond to some resolvent. Lower norms provide a reformulation that
applies both to the global unbounded operator and to those local finite
sections. In Hilbert space, the lower norm of an operator is another name for its smallest singular value:

\begin{definition}[Lower norm]
\label{def:lower-norm}
Let $T\in\cC(\cH_1,\cH_2)$. Its lower norm is defined as
\[
  \nu(T) \coloneqq \inf_{\substack{\psi\in\mathcal D(T)\\ \|\psi\|_{\mathcal H_1}=1}} \|T\psi\|_{\mathcal H_2} = \inf_{\substack{\psi\in\mathcal D(T)\\ \psi \neq 0}} \frac{\|T\psi\|_{\mathcal H_2}}{\|\psi\|_{\mathcal H_1}}.
\]
\end{definition}

The condition $\nu(T)>0$ means that $T$ is bounded below. This implies that
$T$ is injective, and that since $T$ is closed, its range is closed. On the other hand, a closed range and injectivity imply $\nu(T)>0$ by the closed graph theorem. If $T$ is
bijective, then the lower equals the operator norm of the inverse:
\begin{equation}\label{eq:LowerNormEqualsInvOpNorm}
  \nu(T)=\inf_{\substack{\psi\in\mathcal D(T)\\ \|\psi\|_{\mathcal H_1}=1}} \|T\psi\|_{\mathcal H_2}=
  \frac{1}{\sup_{\substack{\psi\in\mathcal D(T)\\ \psi\neq 0}}\frac{\|\psi\|_{\mathcal H_1}}{\|T\psi\|_{\mathcal H_2}}} = \frac{1}{\sup_{\substack{\phi\in \mathcal{H}_2\\ \phi\neq 0}}\frac{\|T^{-1}\phi\|_{\mathcal H_1}}{\|T(T^{-1}\phi)\|_{\mathcal H_2}}} = 
  \|T^{-1}\|^{-1}.
\end{equation}
This can be used to characterise the resolvent set of
$H\in \cC(\cH)$.
However, if $H$ is not self-adjoint, the lower norm of $H-\lambda$
alone does not suffice: it detects injectivity and
closed range, but not surjectivity. This missing information is provided by
the adjoint since $\overline{\operatorname{ran}(H)} = \ker(H^*)^\perp$.

We therefore introduce the combined lower-norm
\[
  \nu_\mathrm{comb}(H) \coloneqq \min\left\{ \nu(H), \nu(H^*) \right\}, \qquad \lambda\in\mathbb C.
\]
We then have $\nu_\mathrm{comb}(H)>0$ if and only if $H^{-1}\in \mathcal B(\mathcal H)$, and $\nu_\mathrm{comb}(H)>0$ implies
$\nu_\mathrm{comb}(H)=\nu(H)$ via~\eqref{eq:LowerNormEqualsInvOpNorm}.
To sum up, we obtain the following alternative characterisation of pseudospectra.
\begin{proposition}[Combined lower norm and resolvent]
\label{prop:combined-lower-norm}
Let $H\in\mathcal C(\mathcal H)$. Then,
\[
  \nu_\mathrm{comb}(H-\lambda) =
  \begin{cases}
  \|(H-\lambda)^{-1}\|^{-1}=\nu(H-\lambda), & \lambda\in\rho(H),\\[1ex] 
    0, & \lambda\in\sigma(H).
  \end{cases}
\]
Consequently, for every $\varepsilon>0$,
\[
  \sigma_\varepsilon(H) = \left\{\lambda\in\mathbb C:\nu_\mathrm{comb}(H-\lambda)<\varepsilon\right\},
\]
and, for every $\varepsilon \geq 0$,
\[
    \Sigma_\varepsilon(H) = \left\{\lambda\in\mathbb C:\nu_\mathrm{comb}(H-\lambda)\leq \varepsilon\right\}.
\]
For $\varepsilon = 0$, we have
\[
    \sigma_0(H)=\Sigma_0(H)=\sigma(H) = \left\{
    \lambda \in \C: \nu_\mathrm{comb}(H-\lambda)=0
    \right\}.
\]
\end{proposition}

For later use, we also recall the following well-known continuity property. 
\begin{proposition}[Continuity of the lower norm] 
\label{prop:cont-lowernorm}
    The lower norm function $\nu\colon \mathcal B(\mathcal H)\to [0,\infty)$ (and consequently $\nu_\mathrm{comb}$) is $1$--Lipschitz continuous.
    For $H\in\mathcal{C}(\mathcal{H})$, the mappings
    \[
  \lambda\longmapsto\nu(H-\lambda),\qquad
  \lambda\longmapsto\nu(H^*-\overline\lambda),\qquad\text{and}\qquad
  \lambda\longmapsto\nu_\mathrm{comb}(H-\lambda)
\]
are $1$--Lipschitz continuous.
\end{proposition}
\begin{proof}
For $\lambda_1,\lambda_2\in\mathbb{C}$ and 
$H\in\mathcal C(\mathcal H)$, we have
\begin{multline*}
\lvert\lambda_1-\lambda_2\rvert=
\lVert Hx-\lambda_2 x-Hx+\lambda_1 x\rVert_{\mathcal{H}}\\
\geq\|(H-\lambda_2)x\|_\mathcal{H}-\|(H-\lambda_1)x\|_\mathcal{H}\geq \nu(H-\lambda_2)-\|(H-\lambda_1)x\|_\mathcal{H}
\end{multline*}
for all $x\in\mathcal{D}(H)$ with $\|x\|_\mathcal{H}=1$, and $\|(H-\lambda_1)x\|_\mathcal{H}$ can approximate
$\nu(H-\lambda_1)$ by definition.
The bounded case works similarly and the minimum of two
$1$--Lipschitz continuous functions is $1$--Lipschitz continuous.
\end{proof}

We will use the last function from Proposition \ref{prop:cont-lowernorm} throughout and refer to it by
\[
\alpha_H \colon\qquad \C \longrightarrow [0,\infty),\qquad \lambda\longmapsto\nu_\mathrm{comb}(H-\lambda)
\]
from here. By the continuity of $\alpha_H$, the pseudospectra $\sigma_\varepsilon(H)$ and $\Sigma_\varepsilon(H)$ are open and closed, respectively, for $\varepsilon >0$.

\subsection{Hausdorff continuity in the bounded case.}
When discussing convergence of set--valued spectral approximations, we need a suitable notion of distance between subsets of $\C$. We motivate our choice by first considering nonempty bounded sets $S,T\subset\C$ and defining\footnote{For a Banach space $X$, $x\in X$ and $Y\subset X$, we define $\dist(x,Y)\coloneqq \inf_{y\in Y} \lVert x-y\rVert_X$.}

\begin{equation}
    \label{eq:HausdorffDefin}
  d_H(S,T)
  \coloneqq \max\Bigl\{
      \sup_{s\in S}\dist(s,T),\,
      \sup_{t\in T}\dist(t,S)
    \Bigr\}
  .
\end{equation}

It is a well-known fact (see, e.g., \cite[Lemma~1.5.1]{Beer1993}) that
\begin{equation}\label{eq:HausdorffEquivDefin}
   d_H(S,T) =
  \sup_{z\in\C}
  \bigl|\dist(z,S)-\dist(z,T)\bigr|.
\end{equation}
This is the Hausdorff distance, a complete metric on the compact subsets of $\C$ and, more generally, a pseudometric on the bounded subsets. In particular, it depends only on the closures of the sets:

$$
  d_H(S,T)=0
  \quad\Longleftrightarrow\quad
  \overline{S}=\overline{T}.
$$

Let now $H\in\cB(\cH)$. Since the ``$\leq$''--pseudospectra
$\Sigma_\varepsilon(H)$ are nonempty compact subsets of $\C$, it is natural to ask whether the map
$[0,\infty)\ni\varepsilon \longmapsto \Sigma_\varepsilon(H) $
is continuous with respect to the Hausdorff metric. The following Lemma shows that continuity from above in the Hausdorff metric follows from $\Sigma_{\varepsilon_0}(H) = \bigcap_{\varepsilon>\varepsilon_0} \Sigma_\varepsilon(H)$ and compactness of the closed pseudospectra. We will provide a proof since it shows methods and approaches that we will apply later in the unbounded case.

\begin{lemma}[Hausdorff limits of monotone compact families]
    Let $(K_n)_{n\in\N}$ and $K$ be nonempty compact subsets of $\C$. Suppose either that $(K_n)$ is decreasing with $K = \bigcap_{n=1}^\infty K_n$, or that $(K_n)$ is increasing with $K=\overline{\bigcup_{n=1}^\infty K_n}$. Then 
    \[
    d_H(K_n,K)\longrightarrow 0
    \]
\end{lemma}
\begin{proof}
    In both cases, one of the two terms in
    \eqref{eq:HausdorffDefin} vanishes, since either
    $K\subset K_n$ or $K_n\subset K$.

    In the decreasing case, suppose that convergence fails. Then, there exists $\delta > 0$ and infinitely many $z_j \in K_{n_j}$ with $\dist(z_j, K) \geq \delta$. By compactness of $K_1$, a subsequence converges to some $z\in K_1\setminus K$. For every fixed $m\in\mathbb{N}$, this subsequence eventually belongs to $K_m$, so $z \in \bigcap_m K_m = K$ yields a contradiction. 

    In the increasing case, compactness of $K$ and density of $\bigcup_n K_n$ in K give, for every $\delta > 0$, a finite $\delta$-net for K in $\bigcup_n K_n$. This net belongs to $K_N$ for some $N\in\mathbb{N}$. Since $K_N \subset K_n \subset K$, for $n \geq N$, we obtain $d_H(K_n,K) \leq d_H(K_N,K) < \delta$.
\end{proof}

For $\varepsilon > 0$, Hausdorff continuity from below follows if we additionally have
\begin{equation*}
    \overline{\sigma_\varepsilon(H)}
    =
    \Sigma_\varepsilon(H)
\end{equation*}
since $\bigcup_{0\leq\eta<\varepsilon}\Sigma_\eta(H)
    = \sigma_\varepsilon(H)$.

In order to address this equality, one has to look at the level sets of the combined lower-norm function $\lambda\longmapsto\nu_{\mathrm{comb}}(H-\lambda I)$.
Its continuity always implies $\overline{\sigma_\varepsilon(H)} \subset \Sigma_\varepsilon(H)$.
For the reverse inclusion, points satisfying $\nu_{\mathrm{comb}}(H-\lambda I)=\varepsilon$
must be approximable by points with $\nu_{\mathrm{comb}}(H-\lambda I)<\varepsilon$, meaning that the level set of $\lambda\longmapsto\nu_{\mathrm{comb}}(H-\lambda I)$ for the level $\varepsilon$ must have an empty interior.
By the maximum modulus theorem, this is already sufficient.
In other words, if the resolvent norm is not constant on any non-empty open subset of
$\rho(H)$, then 
\begin{equation}
\label{eq:open-closed-agree}
 \forall \varepsilon \geq 0: \Sigma_\varepsilon(H)
  =
  \overline{\sigma_\varepsilon(H)}.
\end{equation}

For bounded operators, Globevnik \cite{Globevnik1976} showed that the resolvent norm of $H\in\mathcal{B}(\mathcal{H})$ cannot be locally constant. In fact, if we extend the definition of resolvent and pseudospectrum to the Banach space setting and if $H$ is a bounded operator on a complex Banach space, where either the Banach space or its dual is complex
uniformly convex, then the same conclusion holds. In particular, \eqref{eq:open-closed-agree} and the resulting
Hausdorff continuity hold for every bounded operator on any Hilbert space; see
\cite{Globevnik1976,Shargorodsky2008,Davies_2015}. Hence, we have the following well-known result.

\begin{proposition}[Continuity of pseudospectra in $\varepsilon$]
    Let $H$ be a bounded operator on a complex Banach space $\mathcal{X}$. Then $[0,\infty)\ni\varepsilon \longmapsto \Sigma_\varepsilon(H)$ is continuous in $\varepsilon_0$, meaning 
    \[
      \varepsilon_n\longrightarrow\varepsilon_0 \geq 0
      \quad\Longrightarrow\quad
      d_{\mathrm H}\bigl(
        \Sigma_{\varepsilon_n}(H),
        \Sigma_{\varepsilon_0}(H)
      \bigr)
      \longrightarrow0,
    \]
    whenever $\overline{\sigma_{\varepsilon_0}(H)}=\Sigma_{\varepsilon_0}(H)$. In particular, continuity holds for all $\varepsilon \geq 0$ if $\mathcal{X}$ or its dual is complex uniformly convex, e.g., if $\mathcal{X}$ is a Hilbert space.  
\end{proposition}

This conclusion is false on general Banach spaces. Shargorodsky constructed
a bounded operator whose resolvent norm is constant on a non-empty open set
\cite{Shargorodsky2008}. Such an example may even be realised on a separable,
reflexive, and strictly convex Banach space
\cite{ShkarinShargorodsky2009}. At the corresponding pseudospectral level,
the inclusion
\[
  \overline{\sigma_\varepsilon(H)}
  \subset
  \Sigma_\varepsilon(H)
\]
is strict, and Hausdorff continuity in $\varepsilon$ may fail.

The Lipschitz continuity of $\nu\colon \mathcal B(\mathcal H)\to [0,\infty)$ (Proposition~\ref{prop:cont-lowernorm}) is the main ingredient of the following well-known result.
\begin{proposition}[{Continuous dependence of pseudospectra on operators}]
    Consider $H\in \mathcal B(\cH)$ with $\lambda\mapsto \lVert (H-\lambda)^{-1}\rVert$ not constant on any open subset of $\rho (H)$, and a sequence $(H_n)_{n\in\N}$ from $\mathcal B (\cH)$ that converges to $H$ in the operator norm.
    Then,
    \[
    \lim_{n\to\infty}d_{\mathrm{H}}(\Sigma_\varepsilon(H_n),\Sigma_\varepsilon(H))=0
    \]
    for all $\varepsilon>0$.

\end{proposition}

\subsection{Unbounded operators.}
We assume
that $\rho(H)\neq\varnothing$. Although locally constant
resolvent norms are excluded for bounded operators on Hilbert spaces, this
is no longer true for arbitrary unbounded operators. The counterexample is a block-diagonal closed, densely defined operator
on $\ell^2(\N)$, with $2\times2$ blocks, whose resolvent norm is constant
in a neighbourhood of the origin; see \cite{Shargorodsky2008}. Thus, even
in the Hilbert space setting, one may have
\[
  \overline{\sigma_\varepsilon(H)}
  \subsetneq
  \Sigma_\varepsilon(H)
\]
at a suitable pseudospectral level.

There is, however, a restriction on the possible plateau values. Set
\[
  m_H \coloneqq \inf_{\lambda\in\rho(H)} \|(H-\lambda)^{-1}\|.
\] 
By \cite[Thm.~3.2]{B_gli_2014}, on a complex uniformly convex Banach space, any value attained constantly by the resolvent norm on a non-empty open subset of
$\rho(H)$ must be its global minimum $m_H$. This extends the corresponding results of Globevnik and Vidav \cite[Lemmas~1.1 and~3.0]{GlobevnikVidav1974} and applies, in particular, to Hilbert spaces. 
An open plateau, if it exists, must therefore have height $m_H >0$ since the resolvent norm is strictly positive on $\rho(H)$. Set
\begin{equation} \label{eq:eps*}
  \varepsilon_*\coloneqq\frac1{m_H}
\end{equation}
(possibly $\infty$) as the \emph{critical level}.
For
$\varepsilon>\varepsilon_*$, we infer 
 \[\sigma_\varepsilon(H)=\Sigma_\varepsilon(H)=\C,\]
 and, for $0\leq\varepsilon<\varepsilon_*$,
  \[\overline{\sigma_\varepsilon(H)}=\Sigma_\varepsilon(H)\subsetneq\C.\]
Thus, as $\varepsilon$ decreases, the pseudospectra change from the whole plane to a proper subset at the critical level $\varepsilon_*$, and this is the only possible plateau, i.e.,
\[
    \overline{\sigma_{\varepsilon_*}(H)}
    \subset
    \Sigma_{\varepsilon_*}(H)=\C
\]
is possibly a strict inclusion.
Overall, the critical level, if it exists, does not affect the spectral approximation as $\varepsilon \downarrow 0$.

\begin{convention}
Throughout this paper, we exclude the cases $\sigma(H)=\varnothing$ and $\sigma(H)=\mathbb C$.
Thus, both the spectrum and the resolvent set are assumed to be non-empty. This ensures that the spectral and pseudospectral approximation problems considered below are non-trivial, and it avoids repeatedly separating these exceptional cases.

Moreover, the preceding observations show that the resolvent norm can admit at most one plateau. In limiting arguments with \(\varepsilon\downarrow0\), this exceptional level is eventually passed and therefore becomes irrelevant. We shall hence assume in such arguments that \(\varepsilon>0\) is already sufficiently small to lie below the plateau level, and we will not repeat this condition in each individual theorem.\hfill$\blacksquare$
\end{convention}

It is therefore natural to ask if, and in which sense the pseudospectra
$\sigma_\varepsilon(H),\,  \Sigma_\varepsilon(H)$ converge to $\sigma(H)$ as
$\varepsilon\downarrow0$. Since these sets may be unbounded, the usual
Hausdorff metric is not immediately available, and an appropriate notion of
convergence for unbounded closed subsets of $\C$ has to be chosen. Hausdorff convergence amounts to uniform convergence of the associated distance functions on $\C$. For unbounded sets, this global requirement is generally too restrictive, and a natural replacement would be to require
uniform convergence of the distance functions only on bounded subsets of
$\C$. This leads to the Attouch--Wets metric.

\begin{definition}[Attouch--Wets convergence]
\label{def: Attouch--Wets}
Let \(A,B\subset\C\) be non-empty sets. Their Attouch--Wets distance (cf.~\cite{Mosco1969,AttouchWets1983,AttouchWets1986})
is defined by

\[
d_{\mathrm{AW}}(A,B)
:=
\sum_{k=1}^{\infty}
2^{-k}
\min\left\{
1,\,
\sup_{\lambda\in \overline{B_k(0)}}
\bigl|
\dist(\lambda,A)-\dist(\lambda,B)
\bigr|
\right\}.
\]
By \cite[Lemma~3.1.1]{Beer1993}, this defines a complete metric on the non-empty closed subsets of $\C$.
Since $\dist(\lambda,A)=\dist(\lambda,\overline A)$, the Attouch--Wets distance does not distinguish a set from its closure. Thus, it is only a pseudometric on the non-empty subsets of $\C$.
\end{definition}

\begin{lemma}
    \label{lem:aw-characterization}
For non-empty closed sets \(A_n,A\subset\C\), the following are equivalent:
\[
\begin{aligned}
&\text{\rm (i)} && d_{\mathrm{AW}}(A_n,A)\longrightarrow0,\\
&\text{\rm (ii)} &&
\dist(\,\cdot\,,A_n)\longrightarrow\dist(\,\cdot\,,A)
\quad\text{locally uniformly on }\C,\\
&\text{\rm (iii)} &&
\dist(\lambda,A_n)\longrightarrow\dist(\lambda,A)
\quad\text{for every }\lambda\in\C,\\
&\text{\rm (iv)} && \text{the statement \eqref{eqAWequivalentBoegliConv} below holds:}
\end{aligned}
\]
\begin{multline}\label{eqAWequivalentBoegliConv}
\forall K\subset \C\text{ compact with } \overline{\interior(K)\cap A}=K\cap A\neq\varnothing:\\
\text{ eventually } K\cap A_n\neq \varnothing\text{ and } d_H(A_n\cap K,A\cap K)\longrightarrow 0\text{.}
\end{multline}
\end{lemma}
\begin{proof}
The equivalence of \emph{(i)} and \emph{(ii)} is immediate from the definition, and \emph{(ii)} clearly implies \emph{(iii)}. Conversely, each distance function $\dist(\cdot,A_n)$ is $1$-Lipschitz. Hence, by an Arzelà--Ascoli type argument, pointwise convergence together with equicontinuity implies uniform convergence on compact subsets of $\C$.

Next, assume~\eqref{eqAWequivalentBoegliConv}. We first show that only countably many $\overline{B_{\delta}(0)}$ are not admissible to~\eqref{eqAWequivalentBoegliConv} in place of $K$. Hence, for any $k\in\N$, there exists a radius $\varepsilon_k\in\R$ with $\varepsilon_k\geq k$ such that $\overline{B_{\varepsilon_k}(0)}$ is an admissible $K$:
For $\delta> 0$,
\[z\in\left(\overline{B_{\delta}(0)}\cap A\right)\setminus\overline{B_{\delta}(0)\cap A}\]
readily implies $\lvert z\rvert=\delta$ and the existence of an open ball $z\in B_q(r+\mathrm{i}s)$ with $q,r,s\in\mathbb{Q}$ such that
\[
B_q(r+\mathrm{i}s) \cap B_\delta(0)\cap A = \varnothing\text{,}
\]
which implies
\[
z \in A\cap B_q(r+\mathrm{i}s)\subset \left(B_{\delta}(0)\right)^\complement\text{.}
\]
Therefore, $\delta=\min \{|y|: y \in B_q(r+\mathrm{i}s)\cap A \}$, and the ball $B_q(r+\mathrm{i}s)$ cannot correspond to another radius $\delta^\prime\neq\delta$ in a similar way, which shows the claim about countably many non-admissable $K=\overline{B_{\varepsilon_k}(0)}$. 

For $\lambda\in\C$ and $k\in\N$,~\eqref{eq:HausdorffEquivDefin} now yields
\[
\dist\left(\lambda,A_n\cap \overline{B_{\varepsilon_{k}}}(0)\right)\longrightarrow\dist\left(\lambda, A\cap \overline{B_{\varepsilon_{k}}}(0)\right)\text{.}
\]
Since $A$ is closed and a subset of $\C$, we can find $a\in A$ with
$\lvert\lambda-a\rvert=\dist(\lambda,A)$, i.e., $\dist(\lambda, A\cap \overline{B_{\varepsilon_{k}}}(0))=\dist(\lambda,A)$ for $k\in\N$ with $\overline{B_{\dist(\lambda,A)}}(\lambda)\subset \overline{B_{\varepsilon_{k}}}(0)$.
For a given $\delta>0$, choose $k\in\N$ large enough such that 
\[ 
\dist(\lambda,A_n)\leq\dist(\lambda,A_n\cap\overline{B_{\varepsilon_k}(0)})\leq \dist(\lambda,A)+\delta.
\]
So we have $\dist(\lambda,A_n)=\dist(\lambda,A_n\cap\overline{B_{\varepsilon_k}(0)})$ for $k \in \N$ with 
$\overline{B_{\dist(\lambda,A)+\delta}}(\lambda)\subset\overline{B_{\varepsilon_{k}}}(0)$. For $n\in\N$ sufficiently large,
we obtain
\[
\delta>\left\lvert\dist(\lambda,A\cap\overline{B_{\varepsilon_k}(0)})-\dist\left(\lambda,A_n\cap \overline{B_{\varepsilon_{k}}(0)}\right)\right\rvert=\lvert\dist(\lambda,A)-\dist(\lambda,A_n) \rvert\text{,}
\]
which implies~\emph{(iii)}.

Conversely, assume \emph{(iii)} and let $K\subset \C$ be compact with $\overline{\interior(K)\cap A}=K\cap A\neq \varnothing$. First, we show that $K\cap A_n\neq \varnothing$ for $n\in\N$ large enough (cf.~\cite[Lemma~2.1.2]{Beer1993}). As $\interior(K)\cap A\neq\varnothing$, some $\lambda\in A$ and $\varepsilon>0$ with $\overline{B_\varepsilon}(\lambda)\subset K$ exist. From $\dist(\lambda,A_n)\longrightarrow\dist(\lambda,A)=0$, we deduce 
$\varnothing\neq \overline{B_\varepsilon}(\lambda)\cap A_n \subset K\cap A_n$ eventually for $n\in\N$. For the Hausdorff-convergence, we use \eqref{eq:HausdorffDefin} and first proceed indirectly. Assuming
\[
\sup_{\lambda\in K\cap A_n}\dist(\lambda,K\cap A)\not\longrightarrow 0\text{,}
\]
we find a subsequence $(A_{n_k})_{k\in\N}$, and we can construct a convergent sequence $(\lambda_{k})_{k\in\N}$ with $\lambda_k\in K\cap A_{n_k}$ and $\lambda\coloneqq\lim_{k\to\infty}\lambda_k\in K\setminus A$. Due to (compare the Kuratowski-Painlevé convergence \cite[Prop.~5.2.9]{Beer1993}) $\dist(\lambda,A_{n_k})\leq\vert\lambda-\lambda_{n_k}\rvert
\longrightarrow 0$, this contradicts~\emph{(iii)}. 

Secondly, by compactness of $K\cap A$, for a given $\delta>0$, we can find some $N(\delta)\in\N$, $\lambda_1^\prime,\dots,\lambda_{N(\delta)}^\prime\in \interior(K)\cap A$,
and $\varepsilon_j<\delta$ for $j=1,\dots,N(\delta)$ such that
\[
\sup_{\lambda\in K\cap A}\inf \big\lbrace\lvert\lambda-\lambda_j^\prime\rvert : j=1,\dots,N(\delta)\big\rbrace < \delta\quad \text{ and }\quad \bigcup_{j=1}^{N(\delta)}\overline{B_{\varepsilon_j}(\lambda_j^\prime)}\subset K\text{.}
\]
For sufficiently large $n\in\N$,~\emph{(iii)} yields
\[
\sup_{j=1,\dots,N(\delta)}\dist(\lambda^\prime_j,K\cap A_n) < \delta\implies 
\sup_{\lambda\in K\cap A}\dist(\lambda,K\cap A_n)<2\delta\text{,}
\]
which finishes the proof.
\end{proof}
Note that~\emph{(iii)} is exactly convergence in the Wijsman topology (cf.~\cite{Wijsman1966}).
Further, note that~\eqref{eqAWequivalentBoegliConv} is meant to
incorporate \cite[Theorem~5.3]{Hansen2011} and \cite[Theorem~2.1]{B_gli_2014} into our setting. This requires some additional discussion:
\begin{remark}
\label{rmk:aw-characterization}
    Consider the context of Lemma~\ref{lem:aw-characterization}, i.e., $A_n,A\subset\C$ are always closed and non-empty.
    \begin{itemize}
        \item In fact, we have proven that~\eqref{eqAWequivalentBoegliConv} (and therefore Attouch-Wets convergence) is equivalent to
        \begin{multline}\label{eqAWequivalentBoegliConvOnBalls}
\forall k\in\N\,\exists\,\varepsilon_k\geq k \text{ with } \overline{B_{\varepsilon_k}(0)}\cap A\neq \varnothing:\\
\text{ eventually } \overline{B_{\varepsilon_k}(0)}\cap A_n\neq \varnothing\text{ and } d_H(A_n\cap \overline{B_{\varepsilon_k}(0)},A\cap \overline{B_{\varepsilon_k}(0)})\longrightarrow 0\text{.}
\end{multline}
        \item Note that~\eqref{eqAWequivalentBoegliConv} is equivalent to but different from asking for (see, e.g., \cite[Theorem~3.1.7]{Beer1993} or the recent \cite[Lemma~1.2.4]{Colbrook2026})
        \[
        \max\Bigl\{
      \sup_{a\in A\cap K}\dist(a,A_n),\,
      \sup_{b\in A_n\cap K}\dist(b,A),\, 0
    \Bigr\}\longrightarrow 0
        \]
        for every compact $K\subset\C$. In particular, this shows that, for bounded sets $A,A_n$ with a uniform bound on the sequence, \emph{(i)-(iv)} are equivalent to Hausdorff convergence $d_H(A_n,A)\longrightarrow 0$.
        \item Let $A=\overline{\interior(A)}$. If $K\subset \C$ is compact with $\overline{\interior(K)\cap A}=K\cap A\neq \varnothing$, then also
        $\overline{K\cap \interior(A)}=K\cap A\neq \varnothing$. If $\overline{\interior(K)}=K\subset\C$ is compact with
        $\overline{K\cap \interior(A)}=K\cap A\neq \varnothing$, then also
        $\overline{\interior(K)\cap A}=K\cap A\neq \varnothing$.
        \item Based on the previous observation and the corresponding proof, we see that, for $A=\overline{\interior(A)}$, 
        \begin{multline}\label{eq:BoegliConvGeneralA}
        \forall K\subset \C\text{ compact with } \overline{K\cap \interior(A)}=K\cap A\neq\varnothing:\\
        \text{ eventually } K\cap A_n\neq \varnothing\text{ and } d_H(A_n\cap K,A\cap K)\longrightarrow 0
        \end{multline}
        implies $d_{\mathrm{AW}}(A_n,A)\longrightarrow 0$. The assumptions for $A$ and $K$ now cover \cite[Theorem~2.1]{B_gli_2014}. 
        \item Even if $A=\overline{\interior(A)}$ and $A_n=\overline{\interior(A_n)}$ for $n\in\N$ are closed pseudospectra of operators with nowhere constant resolvent norm, we cannot
        deduce~\eqref{eq:BoegliConvGeneralA} from $d_{\mathrm{AW}}(A_n,A)\longrightarrow 0$: Set $A\coloneqq \C$, $A_n\coloneqq (B_{1/n}(0))^\complement$ and $K\coloneqq\{0\}$. Then, $d_{\mathrm{AW}}(A_n,A)\longrightarrow 0$ and $K\cap \interior(A)=K\cap A=K$, but
        $K\cap A_n=\varnothing$ for all $n\in\N$.
    \end{itemize}
\end{remark}
Finally, we need to introduce a notion of convergence for unbounded operators. Following \cite[Section~IV.2.4]{Kato1995}, we define
\[
S_H\coloneqq\bigl\{(s,Hs):s\in\mathcal{D}(H), \Vert s\rVert_{\cH_1}^2+\lVert Hs\rVert_{\cH_2}^2=1\bigr\}\subset \cH_1\times\cH_2
\]
for $H\in\cC(\cH_1,\cH_2)$. Note that $S_H$ is always closed, bounded, and non-empty.
For $H_1,H_2\in \cC(\cH_1,\cH_2)$, we define the distance
\[
\hat{d}(H_1,H_2)\coloneqq \max\Bigl\{\sup_{u\in S_{H_1}}\dist(u,S_{H_2}),
\sup_{v\in S_{H_2}}\dist(v,S_{H_1})\Bigr\}
\]
which renders $\cC(\cH_1,\cH_2)$ a metric space. The corresponding convergence is called
\emph{convergence in the generalized sense}.
\begin{lemma}[{\cite[Theorems~2.23 and~2.25]{Kato1995}}]
Consider $H\in \cC(\cH)$. A sequence $(H_n)_{n\in\N}$ from
$\cC(\cH)$ converges to $H$ in the generalized sense if and only if there exists some $\lambda\in\rho(H)$ such that eventually $\lambda\in\rho(H_n)$ and 
\begin{equation}\label{eq:GeneralizedSenseViaResolvents}
\bigr\rVert(H_n-\lambda)^{-1}-(H-\lambda)^{-1}\bigl\lVert\longrightarrow 0\text{.}
\end{equation}
In that case, every $\lambda\in\rho(H)$ is eventually in $\rho(H_n)$, and~\eqref{eq:GeneralizedSenseViaResolvents} holds for all $\lambda\in\rho(H)$.

If $H\in \mathcal B(\mathcal H)$, then $\hat{d}(H_n,H)\longrightarrow 0$ is equivalent to $H_n\in \mathcal B(\mathcal H)$ eventually and $\lVert H_n-H\rVert\longrightarrow 0$.
\end{lemma}

We are now able to discuss the continuous dependence of pseudospectra on both the operators and $\varepsilon$, as well as the geometry of pseudospectra.

\begin{proposition}[Continuity of pseudospectra in $\varepsilon$]
\label{prop:aw-continuity-pseudospectra}
Let $\varepsilon \geq 0 \ \forall n \in \N$. Assume that $H\in \cC(\cH)$ and, if $\varepsilon > 0$, it holds\footnote{This condition is satisfied for all but at most one \(\varepsilon>0\), see \eqref{eq:eps*}.}
\begin{equation}
    \label{eq: pseudospectra equal}
    \overline{\sigma_\varepsilon(H)} \; = \;
    \Sigma_\varepsilon(H).
\end{equation}
Then, for every non-negative sequence $\varepsilon_n\to\varepsilon$,
    \[
        \lim_{n\to\infty}d_{\mathrm{AW}} \bigl( \Sigma_{\varepsilon_n}(H), \Sigma_\varepsilon(H) \bigr) = 0.
    \]
\end{proposition}

\begin{proof}
For brevity, set $A_n\coloneqq\Sigma_{\varepsilon_n}(H)\text{ and } A\coloneqq\Sigma_\varepsilon(H)$.
By Lemma \ref{lem:aw-characterization} we need to show that
\begin{eqnarray}
    \label{eq: prop proof}
    \dist(\lambda,A_n) \longrightarrow \dist(\lambda,A) \qquad
  \text{for every }\lambda\in\C,\\
    \nonumber
    \text{i.e., that}\qquad
    \limsup_{n\to\infty}\dist(\lambda,A_n) \leq \dist(\lambda,A) \leq \liminf_{n\to\infty}\dist(\lambda,A_n)
    \qquad \text{for every }\lambda\in\C.
\end{eqnarray}
\begin{enumerate}[label=Step \arabic*: ]
    \item Fix $\lambda\in\C$. We begin with the $\limsup$. If $\varepsilon_n \downarrow \varepsilon$ (in particular, when $\varepsilon=0$) then the inequality $\limsup_{n\to\infty}\dist(\lambda,A_n) \leq \dist(\lambda,A)$ follows directly from $ A_n \supset A$ for $n\in\N$. For the case
    $\varepsilon_n \uparrow \varepsilon$, let $\delta>0$ and
choose $x\in A$ such that
\[
  |\lambda-x| \leq \dist(\lambda,A)+\frac{\delta}{2}.
\]
By \eqref{eq: pseudospectra equal}, there exists $y\in\C$ such that $|x-y|<\frac{\delta}{2}$ and $\nu_\mathrm{comb}(H-y)<\varepsilon$.

Since $\varepsilon_n\to\varepsilon$, we have $\nu_\mathrm{comb}(H-y)<\varepsilon_n$ for all sufficiently large $n$. Hence, $y\in A_n$ for all sufficiently large $n$, and therefore
\[
  \dist(\lambda,A_n) \leq |\lambda-y| \leq
  |\lambda-x|+|x-y| < \dist(\lambda,A)+\delta.
\]
Since $\delta>0$ was arbitrary, it follows that

\begin{equation}
    \label{eq: prop proof 2}
    \limsup_{n\to\infty}\dist(\lambda,A_n) \leq \dist(\lambda,A).
\end{equation}
\item 

We next prove the inequality about the $\liminf$. Passing to a subsequence, if necessary,
we may assume that
\[
  \dist(\lambda,A_n) \longrightarrow
  \liminf_{j\to\infty}\dist(\lambda,A_j).
\]
For each $n$, choose $x_n\in A_n$ such that $|\lambda-x_n| \leq \dist(\lambda,A_n)+\frac1n$.

By (\ref{eq: prop proof 2}), the sequence
$(\dist(\lambda,A_n))_{n\in\N}$ is bounded from above. Consequently, $(x_n)_{n\in\N}$ is
bounded. After passing to a further subsequence, we may therefore assume
that $x_n\to x$
for some $x\in\C$. Since $x_n\in A_n$, we have $\nu_\mathrm{comb}(H-x_n)\leq\varepsilon_n$.
The continuity of $\alpha_H$ and the convergence
$\varepsilon_n\to\varepsilon$ imply $\nu_\mathrm{comb}(H-x) \leq \varepsilon$,
so that $x\in A$. It follows that
\begin{equation}
\label{eq: prop proof 3}
    \dist(\lambda,A) \leq |\lambda-x| =
  \lim_{n\to\infty}|\lambda-x_n| =
  \liminf_{n\to\infty}\dist(\lambda,A_n).
\end{equation}
Combining
\eqref{eq: prop proof 2} and
\eqref{eq: prop proof 3} proves
\eqref{eq: prop proof}.
\end{enumerate}
\end{proof}

The preceding convergence result does not require the pseudospectra to remain uniformly bounded. This reflects a characteristic feature of the Attouch--Wets topology: it only captures the behaviour of sets restricted to bounded subsets of $\C$, so that additional parts may escape to infinity without affecting convergence. Every bounded connected component of a pseudospectrum intersects the spectrum \cite[Thm.~4.3]{TrefethenEmbree2005}; hence, any component which is disjoint from the spectrum must be unbounded. The following result shows that such spectral-free components cannot persist in any bounded region as $\varepsilon\downarrow0$, but instead necessarily escape to infinity.

\begin{theorem}[Escape of spectral-free pseudospectral components]
\label{thm:escape-components}
Let $H\in\cC(\cH)$, $\varepsilon_0>0$, and let $\mathcal E$ be a connected component of $\sigma_{\varepsilon_0}(H)$ satisfying $\mathcal E\cap\sigma(H)=\varnothing$. Then:
\begin{enumerate}[label=(\roman*)]
\item $\mathcal E$ is unbounded.

\item For every $0<\varepsilon<\varepsilon_0$, the set $\mathcal E\cap\sigma_\varepsilon(H)$ is a union of connected components of $\sigma_\varepsilon(H)$. Each of these components is disjoint from $\sigma(H)$ and hence unbounded.

\item The sets $\mathcal E\cap\sigma_\varepsilon(H)$ escape to infinity as $\varepsilon\downarrow0$ in the sense that for every compact set $K\subset\C$ there exists $\varepsilon_K>0$ such that
\[
    \mathcal E\cap\sigma_\varepsilon(H)\cap K=\varnothing
\]
for all $0<\varepsilon<\varepsilon_K$.

\end{enumerate}
\end{theorem}
\begin{proof}
Assertion \emph{(i)} follows immediately from \cite[Thm.~4.3]{TrefethenEmbree2005}, since every bounded connected component of a pseudospectrum intersects the spectrum.

For \emph{(ii)}, let $0<\varepsilon<\varepsilon_0$, and let $\mathcal D$ be a connected component of $\sigma_\varepsilon(H)$ such that $\mathcal D\cap\mathcal E\neq\varnothing$. Since $\mathcal D\subset\sigma_\varepsilon(H)\subset\sigma_{\varepsilon_0}(H)$
and $\mathcal D$ is connected, it must be contained in the connected component $\mathcal E$ of $\sigma_{\varepsilon_0}(H)$. Hence $\mathcal E\cap\sigma_\varepsilon(H)$ is a union of connected components of $\sigma_\varepsilon(H)$. Each such component is contained in $\mathcal E$ and is therefore disjoint from $\sigma(H)$; by (i) it is consequently unbounded.

For \emph{(iii)}, suppose the assertion were false for some compact set $K\subset\C$. Then, there exist $\varepsilon_n\downarrow0$ and
$\lambda_n\in\mathcal E\cap\sigma_{\varepsilon_n}(H)\cap K.$
After passing to a subsequence, compactness of $K$ yields $\lambda_n\to\lambda\in K$. Since $\gamma_H(\lambda_n)<\varepsilon_n$ and $\gamma_H$ is continuous, we obtain $\gamma_H(\lambda)=0$, and hence $\lambda\in\sigma(H)$. On the other hand, $\mathcal E$ is closed relative to $\sigma_{\varepsilon_0}(H)$. Since $\lambda_n\in\mathcal E$ and $\lambda\in\sigma(H)\subset\sigma_{\varepsilon_0}(H),$
it follows that $\lambda\in\mathcal E$, contradicting $\mathcal E\cap\sigma(H)=\varnothing$.
\end{proof}

\begin{proposition}
    The lower norm $\nu\colon \cC(\cH_1,\cH_2) \to [0,\infty)$ is locally Lipschitz continuous with respect to $\hat{d}$, meaning that for every $H \in \cC(\cH_1,\cH_2)$ there exists a neighborhood $U_H$ and a constant $L_{H}>0$ such that \[
    |\nu(V)-\nu(W)|\leq L_H\, \hat{d}(V,W) \quad \forall V, W \in U_H.
    \]
\end{proposition}
\begin{proof}
    We first show that $\nu$ is locally bounded. Fix $H$ and choose $(h,Hh)\in S_H$. Denote $a\coloneqq\|h\|$. If $\hat{d}(H,V) < a/2$, there exists $(v,Vv) \in S_V$ such that $\|(h,Hh)-(v,Vv)\| < a/2$.
    Hence, $\|v\|>a/2$ (triangle inequality) and $\|Vv\|\leq 1$ give $$\nu(V) \leq \|Vv\|/\|v\|\leq 2/a \qquad \forall V\in B_{a/2}(H).$$
    Next, put $M\coloneqq2/a$ and $c\coloneqq\frac1{\sqrt{1+(M+1)^2}}$. \newline
    Let $V,W \in B_{\min\{a/2,c/4\}}(H)$, so $\nu(V),\nu(W)\leq M,$ and choose a minimizing sequence $(v_n,Vv_n)\in S_V$ with
    \begin{equation*}
        \frac{\|Vv_n\|}{\|v_n\|}\longrightarrow\nu(V),
        \qquad
        \frac{\|Vv_n\|}{\|v_n\|}\leq M+1.
    \end{equation*}
    Since $\|v_n\|^2+\|Vv_n\|^2=1$, we have
    \begin{equation}\label{eq:boundfirstcomponent}
       \forall n \in \N: \|v_n\|
        =
        \left(1+\frac{\|Vv_n\|^2}{\|v_n\|^2}\right)^{-1/2}
        \geq c. 
    \end{equation}
    
    Set $d\coloneqq\hat d(V,W)$. If $d=0$, the assertion is immediate. Otherwise, the triangle inequality gives $0<d<c/2$ and we can choose an integer $m\geq2$ such that
    \begin{equation}\label{eq:choiceofm}
        \frac{c}{m+1}\leq d<\frac{c}{m}.
    \end{equation}

    For every $n$, there exists $(w_n,Ww_n)\in S_W$ with
    \[
        \|(v_n,Vv_n)-(w_n,Ww_n)\|<\frac{c}{m}.
    \]
    In particular,
    \[
        \|v_n-w_n\|<\frac{\|v_n\|}{m},
        \qquad
        \|Vv_n-Ww_n\|<\frac{\|v_n\|}{m}.
    \]
    Consequently,
    \[
        \|w_n\|\geq\frac{m-1}{m}\|v_n\|,
        \qquad
        \|Ww_n\|\leq\|Vv_n\|+\frac{\|v_n\|}{m},
    \]
    so
    \[
        \nu(W)
        \leq\frac{\|Ww_n\|}{\|w_n\|}
        \leq
        \frac{1}{m-1}
        +\frac{m}{m-1}\frac{\|Vv_n\|}{\|v_n\|}.
    \]
    Passing to the limit yields
    \begin{equation*}
        \nu(W)
        \leq\frac1{m-1}+\frac{m}{m-1}\nu(V).
    \end{equation*}
    Because $\nu(W)\leq M$, a minimizing sequence for $W$ can be chosen with the same
    lower bound $c$ in \eqref{eq:boundfirstcomponent}.
    Interchanging $V$ and $W$ therefore gives
    \[
        |\nu(V)-\nu(W)|
        \leq\frac{1+\max(\nu(V),\nu(W))}{m-1}
        \leq\frac{1+M}{m-1}.
    \]
    By $m\geq 2$ and \eqref{eq:choiceofm}, we have 
    \[
        \frac1{m-1}
        =
        \frac{m+1}{m-1}\frac1{m+1}
        \leq\frac{3d}{c}.
    \]
    Overall, we obtain
    \begin{equation*}
       \forall V,W \in B_{\min\{a/2,c/4\}}(H): |\nu(V)-\nu(W)| \leq\frac{3(1+M)}{c}\,\hat d(V,W).
    \end{equation*}
    
\end{proof}
\begin{example}
    The lower norm $\nu: \cC(\cH) \to [0,\infty)$ is not uniformly continuous, hence not Lipschitz continuous, with respect to $\hat{d}$.
    \begin{proof}
        Let $I$ denote the identity on $\mathcal H$. For $m\geq 2$ and $(s,ms) \in S_{mI}$, set 
        \[t\coloneqq \underbrace{\sqrt{\frac{1+m^2}{1+(m-1)^2}}}_{C_m\coloneqq}\cdot s,\]
        so that $(t,(m-1)t) \in S_{(m-1)I}$. 
        By the definition of $S_{mI}$, we have $\|s\|^2 = \frac{1}{1+m^2}$.
        Since $C_m$ gives rise to a bijection between $S_{mI}$ and $S_{(m-1)I}$ and $C_m\longrightarrow 1$, we can estimate
        \begin{align*}
            (\hat{d}(mI,(m-1)I))^2 & \leq \sup_{(s,ms) \in S_{mI}} \|ms-(m-1)t\|^2+\|s-t\|^2 \\ & = \sup_{(s,ms) \in S_{mI}} \|s\|^2 \cdot \big(m-(m-1) \cdot C_m\big)^2 + \big(1- C_m\big)^2 \|s\|^2 \\ & = \frac{1}{1+m^2} \cdot \Big(\big(m-(m-1) \cdot C_m\big)^2 + \big(1- C_m\big)^2\Big) \xrightarrow{m \to \infty} 0.
        \end{align*}
        On the other hand, \[
        |\nu(mI)-\nu((m-1)I)|=1.
        \]
        By choosing $\varepsilon = 1/2$, there cannot by a $\delta>0$ such that 
        \[
        \forall m\geq 2:\hat{d}(mI,(m-1)I) < \delta \implies |\nu(mI)-\nu((m-1)I)| < \varepsilon.\qedhere
        \]
    \end{proof}
\end{example}
This local Lipschitz continuity of $\nu\colon \mathcal C(\mathcal H)\to [0,\infty)$ is the main ingredient of the following known result.
\begin{proposition}[{Continuous dependence of pseudospectra on operators; \cite[Thm.~5.3]{Hansen2011} and \cite[Theorem~2.1]{B_gli_2014}}]
    Consider $H\in \cC(\cH)$ with $\lambda\mapsto \lVert (H-\lambda)^{-1}\rVert$ not constant on any open subset of $\rho (H)$, and a sequence $(H_n)_{n\in\N}$ from $\cC(\cH)$ that converges to $H$ in the generalized sense.
    Then,
    \[
    \lim_{n\to\infty}d_{\mathrm{AW}}(\Sigma_\varepsilon(H_n),\Sigma_\varepsilon(H))=0
    \]
    for all $\varepsilon>0$.
\end{proposition}

The representation of pseudospectra through $\alpha_H$ reduces their
approximation to the approximation of lower norms. This is particularly
useful for the localisation procedure developed below, since lower norms
remain meaningful for the rectangular finite-dimensional operators obtained
from local restrictions of $H$ and $H^*$.

%% file: Theorems.tex
\section{Geometric Framework}

\begin{convention}[Standing geometric framework]\label{ass:standing}
Throughout, $G$ is a locally compact, second-countable Hausdorff-group with identity $e$. We assume that $G$ is equipped with a proper left-invariant metric $d$ generating its topology and with a left Haar measure $\mu$. Recall that this measure is, among several other properties, finite on compact sets and non-zero on non-empty open sets. We write $B^\leq_L(x)\coloneqq\{y\in G:d(x,y)\leq L\}.$ We denote by $\Gamma \subset G$ a uniformly discrete subset of $G$, meaning that $\exists r_\Gamma>0$ such that $d(\gamma, \gamma')\geq r_\Gamma$ for all distinct $\gamma,\gamma'  \in \Gamma$. \hfill$\blacksquare$
\end{convention}
Our aim is to treat the continuous and discrete settings simultaneously. Rather than stating the localization results separatly on $L^2(G,\mu)$ and on $\ell^2(\Gamma)$, we introduce a space X, which from now on, is understood to be either the ambient group $G$ or a uniformly discrete (and thus countable) subset $\Gamma \subset G$. However, it is important to understand, that in either case the ambient group $G$ carries the underlying structure and provides the geometry, the localisation centres, and the measure with respect to which the later introduced \emph{localisation functions} are integrated, while no additional group structure is required on $\Gamma$.   

We equip \(X\) with the measure
\[
\mu_X \coloneqq
\begin{cases}
\mu, & X=G,\\
\# , & X=\Gamma,
\end{cases}
\]
where \(\#\) denotes counting measure. Thus,
\[
L^2(X,\mu_X)=
\begin{cases}
L^2(G,\mu), & X=G,\\
\ell^2(\Gamma), & X=\Gamma.
\end{cases}
\]

\begin{notation}
Whenever \(X=G\) or \(X=\Gamma\subset G\) is fixed, $\|\cdot\|_{L^2(X)}\coloneqq\|\cdot\|_{L^2(X,\mu_X)}
\text{ and }
\|\cdot\|_{L^2(G)}\coloneqq\|\cdot\|_{L^2(G,\mu)}$.
 If \(X=\Gamma\), the following localization functions \(\chi_{L,x}\), originally defined on \(G\), are implicitly restricted to \(\Gamma\). Accordingly, multiplication by \(\chi_{L,x}\) and the later appearing commutator \([H,\chi_{L,x}]\) are always understood as operators on \(L^2(X,\mu_X)\). 
\end{notation}

As indicated in point (ii) of the introduction, instead of restricting ourselves to projections $P$, which act as multiplication by characteristic functions, we pass to multiplications of more general localization functions. These localization functions retain compact support and can be understood as "weighted restrictions". 

\begin{definition}[Localization functions]
\label{def:localization-functions}
A measurable function \(\eta\colon [0,\infty)\to[0,1]\) will be called a \emph{localization profile} if \(\eta(0)=1\), \(\eta(t)=0\) for \(t\geq1\), and there exists $\delta \in (0,1]$ such that 
\[
\eta(t)>0 \quad (0\leq t<\delta).
\]
Given such a profile and $L>0$, we define the \emph{family of localization functions at scale $L$} by
\[
\chi_{L,x}(\cdot)\coloneqq\eta\left(\frac{d(x,\cdot)}{L}\right), \qquad x\in G.
\]
The collection $(\chi_{L,x})_{L>0,\,x\in G}$ will be referred to as the \emph{family of localization functions} associated with $\eta$.

By construction, $\chi_{L,x}(x)=1$ and $\operatorname{supp}\chi_{L,x}\subset B_L(x)$. Since $d$ is proper and $\mu$ is finite on compact sets, $\chi_{L,x}\in L^2(G,\mu)$. Since $\chi_{L,x}$ is strictly positive on the open ball $\{ y \in G: d(e,y) < \delta L\}$ with positive Haar measure, we have 
\[
\|\chi_{L,e}\|_{L^2(G)}>0 \quad \text{for every }L>0.
\]

If $X=G$, the functions $\chi_{L,x}$ act on $L^2(G,\mu)$ by pointwise multiplication. If $X=\Gamma\subset G$, we use the same notation for their restrictions to $\Gamma$, acting by pointwise multiplication on $\ell^2(\Gamma)$. In both cases, the localization centres $x$ range over the entire ambient group $G$.
\end{definition}

\begin{notation}
Whenever we refer to \emph{localization functions}, we implicitly assume that they are generated by a fixed localization profile $\eta$ as in Definition~\ref{def:localization-functions}. For simplicity, we use the term \emph{family of localization functions} for the full collection $(\chi_{L,x})_{L>0,\,x\in G}$ and do not repeatedly specify the underlying profile. Likewise, the term \emph{localization function} may refer to any member $\chi_{L,x}$ of this family, with the scale $L>0$ and centre $x\in G$ understood from the context.
\end{notation}

\begin{remark}
The localization technique originates from finite-section methods for infinite matrices. There, finite sections are naturally realized by hard cut-offs, that is, multiplication by characteristic functions of bounded regions. Since such sharp truncations produce boundary terms, Lipschitz cut-offs, in particular tent functions, have proven to be more suitable in more recent localization arguments \cite{HegeMoscolariTeufel2026,LindnerSeifertChandlerWildePreprint,Wittig2026}. We keep the notion of localization profile more general here, since smooth cut-offs will be required later for Schrödinger operators on $L^2(\R^n)$; see Section 6.
\end{remark}
As in the works cited above, the following averaging identity remains a basic ingredient of the localization argument.

\begin{lemma}[Averaging identity]\label{lem:mowing-lawn}
Let \(X=G\) or \(X=\Gamma\subset G\) be as above, and let $(\chi_{L,x})_{L>0,\,x\in G}$ be a family of localization functions. Then, for every \(\psi\in L^2(X,\mu_X)\),
\[
\int_G \|\chi_{L,x}\psi\|_{L^2(X)}^2\,\mathrm d\mu(x)
=
\|\chi_{L,e}\|_{L^2(G)}^2
\|\psi\|_{L^2(X)}^2.
\]
\end{lemma}
\begin{proof}
By Tonelli's theorem,
\[
\int_G \|\chi_{L,x}\psi\|_{L^2(X)}^2\,\mathrm d\mu(x)
=
\int_X |\psi(y)|^2
\left(\int_G |\chi_{L,x}(y)|^2\,\mathrm d\mu(x)\right)
\mathrm d\mu_X(y).
\]
For fixed $y\in X\subset G$, the substitution $x=yz$, together with left invariance of $d$ and $\mu$, gives
\[
\int_G |\chi_{L,x}(y)|^2\,\mathrm d\mu(x)
=
\int_G \left|\eta\left(\frac{d(yz,y)}{L}\right)\right|^2\,\mathrm d\mu(z)
=
\int_G \left|\eta\left(\frac{d(z,e)}{L}\right)\right|^2\,\mathrm d\mu(z)
=
\|\chi_{L,e}\|_{L^2(G)}^2.
\]
Hence,
\[
\int_G \|\chi_{L,x}\psi\|_{L^2(X)}^2\,\mathrm d\mu(x)
=
\|\chi_{L,e}\|_{L^2(G)}^2
\int_X |\psi(y)|^2\,\mathrm d\mu_X(y)
=
\|\chi_{L,e}\|_{L^2(G)}^2\|\psi\|_{L^2(X)}^2.
\]
\end{proof}

\section{Local sections and graph-norm localisation}

The localization profiles allow us to cut off pseudoeigenvectors and
thereby obtain vectors supported in bounded regions. From the operator
point of view, this amounts to restricting \(H\) to vectors in its
domain whose support lies in a prescribed window. For a matrix operator
on \(\ell^2(\Gamma)\), one may think of this as retaining the columns
indexed by \(\Gamma\cap B^\leq_L(x)\): the input vector has finite support,
while the output may have entries anywhere in \(\Gamma\). We therefore
keep the full output space in the general definition.

\begin{definition}[Local sections]\label{def:local-sections}
Let \(X=G\) or \(X=\Gamma\subset G\) be as above, and let
\[
H\colon D(H)\subset L^2(X,\mu_X)\longrightarrow L^2(X,\mu_X)
\]
be closed and densely defined. For \(L>0\) and \(x\in G\), set
\(X_{L,x}\coloneqq X\cap B^\leq_L(x)\).\footnote{If \(X=\Gamma\), this window
may be empty. Then \(L^2(X_{L,x},\mu_X)=\{0\}\), \(H_{L,x}\) is
the unique linear operator from \(\{0\}\) into \(L^2(X,\mu_X)\)
and \(\nu(H_{L,x}-\lambda)=\infty\). Later, such windows won't affect the infimum over \(x\in G\).} Denote by
\[
P_{L,x}\colon L^2(X,\mu_X)\longrightarrow L^2(X_{L,x},\mu_X)
\]
the restriction operator. Its adjoint \(P_{L,x}^*\) is extension by zero.

The \emph{local section} of \(H\) at \(x\) and scale \(L\) is
\[
\begin{aligned}
H_{L,x}\colon D(H_{L,x})&\subset L^2(X_{L,x},\mu_X)
                         \longrightarrow L^2(X,\mu_X),\\
D(H_{L,x})
&\coloneqq\bigl\{u\in L^2(X_{L,x},\mu_X):
          P_{L,x}^*u\in D(H)\bigr\},\\
H_{L,x}u&\coloneqq HP_{L,x}^*u.
\end{aligned}
\]
For \(\lambda\in\mathbb C\), we write
\[
(H_{L,x}-\lambda)u\coloneqq (H-\lambda)P_{L,x}^*u,
\qquad u\in D(H_{L,x}).
\]
\end{definition}

The localization functions connect the global operator with its local
sections. Suppose that multiplication by \(\chi_{L,x}\) preserves
\(D(H)\). For \(\psi\in D(H)\), the vector \(\chi_{L,x}\psi\) belongs
to \(D(H)\) and vanishes outside \(X_{L,x}\). Thus
\(u\coloneqq P_{L,x}(\chi_{L,x}\psi)\) belongs to \(D(H_{L,x})\), and
\[
\|u\|_{L^2(X_{L,x})}=\|\chi_{L,x}\psi\|_{L^2(X)},
\qquad
(H_{L,x}-\lambda)u=(H-\lambda)\chi_{L,x}\psi.
\]
Whenever \(\chi_{L,x}\psi\neq0\), it follows that
\begin{equation}\label{estimateLowerNormLocalSectionLocalizedEigenValue}
\nu(H_{L,x}-\lambda)
\leq
\frac{\|(H-\lambda)\chi_{L,x}\psi\|_{L^2(X)}}
     {\|\chi_{L,x}\psi\|_{L^2(X)}}.
\end{equation}
Recalling the commutator notation, $[A,B]\coloneqq AB-BA$, we now come to the key property for our localization argument:

\begin{definition}[Commutator Property]\label{def:commute}
Let $X=G$ or $X=\Gamma\subset G$ be as above, let $(\chi_{L,x})_{L>0,\,x\in G}$ be a family of localization functions, and let
\[
H\colon D(H)\subset L^2(X,\mu_X)\to L^2(X,\mu_X)
\]
be closed and densely defined. 
Then, we will say that $H$ and the family $(\chi_{L,x})_{L>0,\,x\in G}$ \emph{have the commutator property} if multiplication by $\chi_{L,x}$ preserves $D(H)$ for every $L>0$ and $x\in G$ and if there exist constants\footnote{The threshold $L_0$ allows the commutator estimate to hold eventually. In the Schr\"odinger example of Section~\ref{sec:schrodinger-core-closure}, the Laplacian of the cutoff produces an \(L^{-2}\) term; there we can take \(L_0=1\).} $L_0>0$ and $C_H \geq 0$ such that, for every $L\geq L_0$ and every $\psi\in D(H)$,
\begin{equation} \label{eq:comm-ineq}
\left(\int_G \|[H,\chi_{L,x}]\psi\|_{L^2(X)}^2\,\mathrm d\mu(x)\right)^{1/2}
\leq
\frac{C_H}{L}\|\chi_{L,e}\|_{L^2(G)}
\bigl(\|H\psi\|_{L^2(X)}+\|\psi\|_{L^2(X)}\bigr).
\end{equation}
\end{definition}

In \cite{LindnerSeifertChandlerWildePreprint, Wittig2026} this commutator property is shown to hold for bounded operators $H$ with a finite interaction range or even with a Wiener-type decay under different assumptions on $X$ and the localizing family. Note that, for bounded operators $H$, the graph norm on the right of inequality \eqref{eq:comm-ineq} can be replaced by the $L^2(X)$-norm of $\psi$ at the cost of a different constant $C_H$. In Sections 5 and 6 below, we give conditions on unbounded operators on $\ell^2(\Z)$ as well as on the potential of Schrödinger operators on $L^2(\R^d)$ that guarantee $H$ and the corresponding localizing family to have the commutator property.

If $H$ and $(\chi_{L,x})_{L>0,\,x\in G}$ have the commutator property, we can build the announced bridge between local and global pseudomodes of $H$:

\begin{lemma}
\label{lem:localise-graphnorm}
Let $X=G$ or $X=\Gamma\subset G$ be as above, let $(\chi_{L,x})_{L>0,\,x\in G}$ be a family of localization functions, and let
\[
H\colon D(H)\subset L^2(X,\mu_X)\to L^2(X,\mu_X)
\]
be closed and densely defined. Now, assume that $H$ and $(\chi_{L,x})_{L>0,\,x\in G}$ have the commutator property and let $L_0>0$ and $C_H \geq 0$ denote the corresponding constants from Definition \ref{def:commute}.

Let \(\lambda\in\mathbb C\) and suppose that \(\psi\in D(H)\) satisfies
\[
\|(H-\lambda)\psi\|_{L^2(X)}
\leq
\varepsilon\|\psi\|_{L^2(X)}
\]
for some \(\varepsilon\geq0\). Then, for every $L\geq L_0$, there exists \(x\in G\) with \(\chi_{L,x}\psi\neq0\) such that
\[
\|(H-\lambda)\chi_{L,x}\psi\|_{L^2(X)}
\leq
\left(\varepsilon+\frac{C_H}{L}(\varepsilon+|\lambda|+1)\right)
\|\chi_{L,x}\psi\|_{L^2(X)}.
\]
Consequently, for some \(x\in G\),
\[
\nu(H_{L,x}-\lambda)
\leq
\varepsilon+\frac{C_H}{L}(\varepsilon+|\lambda|+1).
\]
\end{lemma}

\begin{proof}
Fix $L\geq L_0$. For every \(x\in G\),
\[
(H-\lambda)\chi_{L,x}\psi
=
\chi_{L,x}(H-\lambda)\psi+[H,\chi_{L,x}]\psi.
\]
Hence, by Minkowski's inequality, the averaging identity from Lemma~\ref{lem:mowing-lawn}, and the assumed commutator estimate,
\[
\begin{aligned}
\left(\int_G\|(H-\lambda)\chi_{L,x}\psi\|_{L^2(X)}^2\,\mathrm d\mu(x)\right)^{1/2}
&\leq
\left(\int_G\|\chi_{L,x}(H-\lambda)\psi\|_{L^2(X)}^2\,\mathrm d\mu(x)\right)^{1/2} \\
&+
\left(\int_G\|[H,\chi_{L,x}]\psi\|_{L^2(X)}^2\,\mathrm d\mu(x)\right)^{1/2} \\
&\leq
\|\chi_{L,e}\|_{L^2(G)}
\left(
\varepsilon\|\psi\|_{L^2(X)}
+
\frac{C_H}{L}\bigl(\|H\psi\|_{L^2(X)}+\|\psi\|_{L^2(X)}\bigr)
\right).
\end{aligned}
\]
Moreover,
\[
\|H\psi\|_{L^2(X)}
\leq
\|(H-\lambda)\psi\|_{L^2(X)}+|\lambda|\|\psi\|_{L^2(X)}
\leq
(\varepsilon+|\lambda|)\|\psi\|_{L^2(X)},
\]
and therefore
\[
\left(\int_G\|(H-\lambda)\chi_{L,x}\psi\|_{L^2(X)}^2\,\mathrm d\mu(x)\right)^{1/2}
\leq
\left(\varepsilon+\frac{C_H}{L}(\varepsilon+|\lambda|+1)\right)
\|\chi_{L,e}\|_{L^2(G)}\|\psi\|_{L^2(X)}.
\]
Again by Lemma~\ref{lem:mowing-lawn},

\[
\left(\int_G\|(H-\lambda)\chi_{L,x}\psi\|_{L^2(X)}^2\,\mathrm d\mu(x)\right)^{1/2}
\leq
\left(\varepsilon+\frac{C_H}{L}(\varepsilon+|\lambda|+1)\right)
\left(\int_G \|\chi_{L,x}\psi\|_{L^2(X)}^2\,\mathrm d\mu(x)\right)^{1/2}.
\]

Since $\psi\neq0$, $ \int_G \|\chi_{L,x}\psi\|_{L^2(X)}^2\,\mathrm d\mu(x) = \|\chi_{L,x}\|_{L^2(G)}^2\|\psi\|_{L^2(X)}^2 > 0. $
In particular, the set of $x\in X$ for which $\chi_{L,x}\psi\neq0$ has positive measure. If
\[
  \|(H-\lambda)\chi_{L,x}\psi\|_{L^2(X)}
  >
   \left(
        \varepsilon+\frac{C_H}{L}(\varepsilon+|\lambda|+1)
    \right) \|\chi_{L,x}\psi\|_{L^2(X)}
\]
held for every such $x$, then squaring and integrating would contradict the
preceding averaged inequality. Hence there exists $x\in X$ such that $\chi_{L,x}\psi\neq0$
and
\[
  \|(H-\lambda)\chi_{L,x}\psi\|_{L^2(X)}
  \leq
  \Bigl(\varepsilon+\frac{C_H}{L}(\varepsilon+|\lambda|+1) \Bigr)\|\chi_{L,x}\psi\|_{L^2(X)}.
\]
Thus $\chi_{L,x}\psi$ provides a non-trivial localized vector with controlled
residual.

The lower-norm estimate follows from~\eqref{estimateLowerNormLocalSectionLocalizedEigenValue}.
\end{proof}

\begin{corollary}\label{cor:global-local}
Let $X=G$ or $X = \Gamma \subset G$ be as above, let $H:D(H)\subset L^2(X,\mu_X)\to L^2(X,\mu_X)$ be closed and densely defined ,and suppose that the assumptions of Lemma~\ref{lem:localise-graphnorm} are satisfied. Then, for every $\lambda\in\C$ and every $L\geq L_0$,
\begin{equation}
\label{eq:global-local}
    \inf_{x\in G}\nu(H_{L,x}-\lambda) \leq \nu(H-\lambda) \bigl(1+\frac{C_H}{L} \bigr) + \frac{C_H}{L}\bigl(|\lambda|+1\bigr). 
\end{equation} 
\end{corollary}

\begin{proof}
Fix $\lambda\in\C$ and $L>L_0$.  Let $\varepsilon>\nu(H-\lambda)$.  By the definition
of the lower norm, there exists $\psi\in L^2(X)$ such that $\|(H-\lambda)\psi\|_{L^2(X)} < \varepsilon \, \|\psi\|_{L^2(X)}$.
Applying Lemma~\ref{lem:localise-graphnorm} to this vector $\psi$, we obtain some $x\in X$ with $\nu(H_{L,x}-\lambda)\ \le\ \varepsilon+\frac{C_0}{L}$.
Taking the infimum over $x$ gives $\inf_{x\in X}\nu(H_{L,x}-\lambda)\ \le\ \varepsilon+\frac{C_0}{L}$.
Since this estimate holds for every $\varepsilon>\nu(H-\lambda)$, \eqref{eq:global-local} follows by letting $\varepsilon\downarrow\nu(H-\lambda)$.
\end{proof}

\begin{theorem}[Window pseudospectra and convergence]\label{thm:window-pseudo-convergence}
Let $X=G$ or $X=\Gamma\subset G$ be as above, and let $H:D(H)\subset L^2(X,\mu_X)\to L^2(X,\mu_X)$
be closed and densely defined. Suppose that the assumptions of Lemma~\ref{lem:localise-graphnorm} are satisfied for both $H$ and $H^*$, with a common constant $C_H>0$. For $\lambda\in\C$, recall $\nu_\mathrm{comb}(H-\lambda) \coloneqq
\min\bigl\{
\nu(H-\lambda),
\nu(H^*-\overline{\lambda})
\bigr\}$ and define
\[
\alpha_{H,L}(\lambda) \coloneqq
\inf_{x\in G}
\min\bigl\{
\nu(H_{L,x}-\lambda),
\nu((H^*)_{L,x}-\overline{\lambda})
\bigr\}.
\]
Then, for every $L\geq L_0$ and every $\lambda\in\C$,
\begin{equation}\label{eq:two-sided-local-lower-norm}
\alpha_H(\lambda)
\leq
\alpha_{H,L}(\lambda)
\leq
\alpha_H(\lambda)\left(1+\frac{C_H}{L}\right)
+\frac{C_H}{L}\bigl(|\lambda|+1\bigr).
\end{equation}
In particular, $\alpha_{H,L}\to\alpha_H$ locally uniformly on $\C$ as $L\to\infty$. \\
For $\varepsilon>0$ and $\eta\geq0$, let\footnote{In general, the sets $\gamma_{L,\varepsilon}(H)$ and $\Gamma_{L,\varepsilon}(H)$ need not coincide. For a counterexample see \cite[Section~5.2]{ChandlerWildeChonchaiyaLindner2024}.}
\[
\gamma_{L,\varepsilon}(H)\coloneqq
\{\lambda\in\C:\alpha_{H,L}(\lambda)<\varepsilon\},
\qquad
\Gamma_{L,\eta}(H)\coloneqq
\{\lambda\in\C:\alpha_{H,L}(\lambda)\leq\eta\}.
\]
Then the following assertions hold.

\begin{enumerate}
\item[(a)] \textup{(Pseudospectral inclusion bounds)}
For every $\varepsilon>0$ and every $\eta\geq0$,
\[
\gamma_{L,\varepsilon}(H)\subset \sigma_\varepsilon(H),
\qquad
\Gamma_{L,\eta}(H)\subset \Sigma_\eta(H).
\]
Moreover, for every $R>0$,
\[
\sigma_\varepsilon(H)\cap B_R(0)
\subset
\gamma_{L,\varepsilon+\delta_{L,R}(\varepsilon)}(H),
\quad \text{and} \quad
\Sigma_\eta(H)\cap B_R(0)
\subset
\Gamma_{L,\eta+\delta_{L,R}(\eta)}(H),
\]
where
\[
\delta_{L,R}(t)\coloneqq\frac{C_H}{L}(t+R+1).
\]

\item[(b)] \textup{(Attouch--Wets convergence)}
For every fixed $\varepsilon>0$, $\varepsilon \neq \varepsilon_*$ with $\varepsilon_*$ from~\eqref{eq:eps*}, and every fixed $\eta>0$, $\eta \neq \varepsilon_*$,
\[
  \gamma_{L,\varepsilon}(H)\xrightarrow[d_{AW}]{L\to\infty}\sigma_\varepsilon(H),
  \qquad
  \Gamma_{L,\eta}(H)\xrightarrow[d_{AW}]{L\to\infty}\Sigma_\eta(H).
\]
Moreover, for $\widetilde{\Gamma}_{L,0}(H) \coloneqq \left\{ \lambda\in\C: \alpha_{H,L}(\lambda) \leq \frac{C_H}{L}\bigl(|\lambda|+1\bigr) \right\}$, for every $L\geq L_0$,
 
\[
\sigma(H)\subset\widetilde{\Gamma}_{L,0}(H) \quad \text{and} \quad
d_{\mathrm{AW}}
\bigl(
\widetilde{\Gamma}_{L,0}(H),
\sigma(H)
\bigr)
\longrightarrow0
\qquad\text{as }L\to\infty.
\]

\end{enumerate}
\end{theorem}

\begin{proof}
The first inequality in~\eqref{eq:two-sided-local-lower-norm} follows by definition, the second from Corollary~\ref{cor:global-local}.

For \textup{(a)}, the first inclusion follows directly from the lower bound in \eqref{eq:two-sided-local-lower-norm}. Let $\lambda\in\gamma_{L,\varepsilon}(H)$, then $\alpha_{H,L}(\lambda)<\varepsilon$, and hence $\alpha_H(\lambda)\leq\alpha_{H,L}(\lambda)<\varepsilon$.
Thus $\lambda\in\sigma_\varepsilon(H)$, so $\gamma_{L,\varepsilon}(H)\subset\sigma_\varepsilon(H)$.

By the same arguments $\Gamma_{L,\eta}(H)\subset\Sigma_\eta(H)$.
\noindent
For the reverse inclusions, fix $R>0$. Let first
$\lambda\in\sigma_\varepsilon(H)\cap B_R(0)$. Then $\alpha_H(\lambda)<\varepsilon$ and $|\lambda|\leq R$. By the upper bound in \eqref{eq:two-sided-local-lower-norm},
\[
\alpha_{H,L}(\lambda)
\leq
\alpha_H(\lambda)\left(1+\frac{C_H}{L}\right)
+\frac{C_H}{L}\bigl(|\lambda|+1\bigr)
<
\varepsilon+\frac{C_H}{L}(\varepsilon+R+1).
\]
Hence
$\alpha_{H,L}(\lambda)
<
\varepsilon+\delta_{L,R}(\varepsilon)$,
and therefore
$
\lambda\in
\gamma_{L,\varepsilon+\delta_{L,R}(\varepsilon)}(H)
$, so
$\sigma_\varepsilon(H)\cap B_R(0)
\subset
\gamma_{L,\varepsilon+\delta_{L,R}(\varepsilon)}(H)$.

By the same arguments $\Sigma_\eta(H)\cap B_R(0) \subset \Gamma_{L,\eta+\delta_{L,R}(\eta)}(H)$.

For \textup{(b)}, fix $\eta>0$ and set
\[
R_L\coloneqq\sqrt{L},
\qquad
\eta_L^-\coloneqq
\frac{\eta-\frac{C_H}{L}(R_L+1)}
     {1+\frac{C_H}{L}}.
\]
For all sufficiently large $L$, we have $\eta_L^->0$ and
$\eta_L^-\to\eta$. If $\lambda\in\Sigma_{\eta_L^-}(H)$ and
$|\lambda|\leq R_L$, then \eqref{eq:two-sided-local-lower-norm} gives
\[
\alpha_{H,L}(\lambda)
\leq
\left(1+\frac{C_H}{L}\right)\eta_L^-
+\frac{C_H}{L}(R_L+1)
=\eta.
\]
Hence
\begin{equation}
\label{eq: convergence part b}
\Sigma_{\eta_L^-}(H)\cap B_{R_L}(0)
\subset
\Gamma_{L,\eta}(H)
\subset
\Sigma_\eta(H).
\end{equation}

Fix $x\in\C$. By \eqref{eq: convergence part b},
\[
\dist(x,\Sigma_\eta(H))
\leq
\dist(x,\Gamma_{L,\eta}(H))
\leq
\dist\bigl(x,\Sigma_{\eta_L^-}(H)\cap B_{R_L}(0)\bigr).
\]
By Proposition~\ref{prop:aw-continuity-pseudospectra},
$\dist(x,\Sigma_{\eta_L^-}(H))\to\dist(x,\Sigma_\eta(H))$. For each $L$,
choose $y_L\in\Sigma_{\eta_L^-}(H)$ such that
\[
|x-y_L|
\leq
\dist(x,\Sigma_{\eta_L^-}(H))+\frac1L.
\]
The sequence $(y_L)$ is bounded, since the right-hand side is bounded. As
$R_L\to\infty$, we therefore have $y_L\in B_{R_L}(0)$ for all sufficiently
large $L$, and consequently
\[
\dist\bigl(x,\Sigma_{\eta_L^-}(H)\cap B_{R_L}(0)\bigr)
\leq |x-y_L|
\leq \dist(x,\Sigma_{\eta_L^-}(H))+\frac1L
\longrightarrow
\dist(x,\Sigma_\eta(H)).
\]
Thus
$\dist(x,\Gamma_{L,\eta}(H))\to\dist(x,\Sigma_\eta(H))$ for every
$x\in\C$. By the characterization of Attouch--Wets convergence in
Lemma~\ref{lem:aw-characterization}, 
it follows that
\[
d_{\mathrm{AW}}
\bigl(\Gamma_{L,\eta}(H),\Sigma_\eta(H)\bigr)
\longrightarrow0.
\]

Since the Attouch--Wets pseudometric does not distinguish a set from its closure, the same argument applies to the open pseudospectra after replacing $\Sigma_{\eta_L^-}(H)$ by $\sigma_{\eta_L^-}(H)$ and $\Sigma_{\eta_L}(H)$ by $\sigma_{\eta_L}(H)$. Then sandwiching $\sigma_{\eta_L^-}(H) \subset \gamma_{L,\eta}(H) \subset \sigma_{\eta_L}(H)$, yields the claim.

For the last statement, \eqref{eq:two-sided-local-lower-norm} gives 
\[ 
\sigma(H) \subset \widetilde{\Gamma}_{L,0}(H) \subset \left\{ \lambda\in\C: \alpha_{H}(\lambda) \leq \frac{C_H}{L}\bigl(|\lambda|+1\bigr) \right\}.
\]
Let $R>0$ be any radius admissible w.r.t. $A = \sigma$ in \eqref{eqAWequivalentBoegliConvOnBalls}, intersecting with $B_R(0)$ yields
    \[
    \label{eq:Gammetildeinclusion}
    \sigma(H) \cap B_R(0) \subset \widetilde{\Gamma}_{L,0}(H) \cap B_R(0) \subset \sigma_{\frac{C_H}{L}(R+1)}(H) \cap B_R(0)
    \]
Since $d_{AW}(\sigma,\sigma_{\frac{C_H}{L}(R+1)})\to 0$ as $L \to \infty$, Remark \ref{rmk:aw-characterization} implies 
\[
d_H(\sigma_{\frac{C_H}{L}(R+1)}(H) \cap B_R(0), \sigma(H) \cap B_R(0)) \to 0.
\]
A sandwich argument with the set inclusion above and applying the same characterization from Remark \ref{rmk:aw-characterization} once more proves the claim.
\end{proof}

%% file: discrete_example.tex
\section{Discrete Examples}

We now turn to discrete operators $H$ and verify the commutator property from Definition \ref{def:commute} with a family $(\chi_{L,x})$ of so-called tent functions. Let $G$ be as above and assume
additionally that its Haar measure satisfies the doubling
property: there exists $C_D\geq1$ such that
\[
    \mu(B^\leq_{2r}(x))
    \leq C_D\,\mu(B^\leq_r(x)),
    \qquad x\in G,\quad r>0.
\]
Let $\Gamma\subset G$ satisfy the standing geometric
assumptions \ref{ass:standing}, and set $X\coloneqq\Gamma$, equipped with the
restricted metric and counting measure.

We localize using the tent functions
\[
    w_{L,x}(y)\coloneqq
    \max\left\{1-\frac{d(x,y)}{L},0\right\},
    \qquad L>0,\quad x,y\in G.
\]
Their restrictions to $\Gamma$ define the multiplication
operators
\[
    (W_{L,x}\psi)(p)\coloneqq w_{L,x}(p)\psi(p),
    \qquad p\in\Gamma.
\]
These functions are supported in $B^\leq_L(x)$ and $L^{-1}$-Lipschitz. Together with the doubling property,
this provides the geometric control needed for the
averaged commutator estimates below.

Our focus is on possibly unbounded operators of finite
interaction range on $\ell^2(\Gamma)$:
\begin{definition}
    Let $(H_{pq})_{p,q\in\Gamma}$ be a (potentially bi-infinite) matrix in $\C$ with $H_{pq}=0$ whenever $d(p,q)>m$ for some $m\geq 0$.
    Set $v(p)\coloneqq H_{pp}$ and define the \emph{off-diagonal edge energy} by
\[
    \mathcal E_{\mathrm{off}}(\psi)^2
    \coloneqq
    \sum_{p\in\Gamma}
    \sum_{\substack{q\in\Gamma\\0<d(p,q)\leq m}}
    |H_{pq}|^2|\psi(q)|^2\in [0,\infty],
\]
where $\psi\in\ell^2(\Gamma)$.
Further, define the operator $H_{\mathrm{max}}\colon D(H_{\mathrm{max}})\subset \ell^2(\Gamma)\to \ell^2(\Gamma)$
on its maximal $\ell^2$-domain
\[
    D(H_{\mathrm{max}})\coloneqq
    \Bigl\{\psi\in\ell^2(\Gamma): \bigl(\sum_{q\in\Gamma}
    H_{pq}\psi(q)\bigr)_{p\in\Gamma}\in\ell^2(\Gamma) \Bigr\}
\]
via
\[
    (H_{\mathrm{max}}\psi)(p)
    \coloneqq
    v(p)\psi(p)
    +
    \sum_{\substack{q\in\Gamma\\0<d(p,q)\leq m}}
    H_{pq}\psi(q),
    \qquad p\in\Gamma.
\]
    We then call $H_{\mathrm{max}}$ the \emph{maximal operator of finite interaction range} $m$ corresponding to $(H_{pq})_{p,q\in\Gamma}$.
\end{definition}
We first discuss some properties of such an $H_{\mathrm{max}}$.
\begin{remark}\label{rem:PropertiesFiniteInteractionRangeUnbounded}
    Let $H_{\mathrm{max}}$ be the maximal operator of finite interaction range $m\geq 0$ corresponding to a matrix $(H_{pq})_{p,q\in\Gamma}$ in $\C$.
    \begin{itemize}
        \item For the canonical basis $(\delta_i)_{i\in\Gamma}$, we clearly have
        $H_{\mathrm{max}}(\delta_i)\in \ell^2(\Gamma)$. Hence, $c_{00}(\Gamma)\subset D(H_{\mathrm{max}})$, where $c_{00}(\Gamma)$ stands for the
        sequences with finitely many non-zero entries. Define the densely defined and closed operator
        \[
        H_{\mathrm{min}}\coloneqq \overline{H_{\mathrm{max}}\restriction_{c_{00}(\Gamma)}}.
        \]
        One readily sees that $(H^\dagger_{\mathrm{min}})^*=H_{\mathrm{max}}$, where $H^\dagger_{\mathrm{min}}$ corresponds to the adjoint matrix $(\overline{H_{qp}})_{p,q\in\Gamma}$.
        Therefore, $H_{\mathrm{max}}$ is densely defined and closed with
        $H_{\mathrm{min}}\subset H_{\mathrm{max}}$.
        \item For $\psi\in\ell^2(\Gamma)$,
        \[
        \lVert H_{\mathrm{max}}\psi\rVert_{\ell^2(\Gamma)}^2 \leq N_m\Bigl( \sum_{p\in\Gamma}
    \lvert v(p)\rvert^2\lvert\psi(p)\rvert^2+ \sum_{p\in\Gamma}\sum_{\substack{q\in\Gamma\\0<d(p,q)\leq m}}
    |H_{pq}|^2|\psi(q)|^2\Bigr)\in [0,\infty],
        \]
        where $N_m\coloneqq \sup_{p\in\Gamma}\#(B^\leq_m(p)\cap\Gamma)$. By the standing geometric assumptions,
\[
    N_m\mu(B^\leq_{\tau}(e))
    \leq \mu(B^\leq_{m+\tau}(e))\implies 1\leq N_m\leq \frac{\mu(B^\leq_{m+\tau}(e))}{\mu(B^\leq_{\tau}(e))}
    <\infty,
\]
where $\tau>0$ satisfies $B^\leq_\tau(p)\cap \Gamma =\{p\}$ for all $p\in\Gamma$.
Thus, we can define $H\subseteq H_{\mathrm{max}}$ on the domain
\begin{equation}\label{eq:VerySmallDomainFiniteInterRange}
D(H)\coloneqq \left\{\psi\in\ell^2(\Gamma): (v(p)\psi(p))_{p\in\Gamma}\in \ell^2(\Gamma)\text{ and }  \mathcal E_{\mathrm{off}}(\psi)<\infty\right\},
\end{equation}
and $c_{00}(\Gamma)\subset D(H)$ is a core.
Overall, we have
\[
H\subset \overline{H}= H_{\mathrm{min}}\subset H_{\mathrm{max}}.
\]
\item The closed graph theorem implies that
\[
H= \overline{H}= H_{\mathrm{min}}= H_{\mathrm{max}}
\]
if and only if
there exists an $C>0$ such that
\begin{equation}\label{eq:HequalsHmaxequalsHmin}
    \lVert v\psi\rVert_{\ell^2(\Gamma)}^2+ \mathcal E_{\mathrm{off}}(\psi)^2\leq C\big(\lVert\psi\rVert_{\ell^2(\Gamma)}^2+\lVert H_{\mathrm{max}}\psi\rVert_{\ell^2(\Gamma)}^2\big)
\end{equation}
for all $\psi\in D(H_{\mathrm{max}})$. In particular,
the graph norm of $H_{\mathrm{max}}$ is then equivalent to
$\lVert\cdot\rVert_{\ell^2(\Gamma)}+\lVert v\cdot\rVert_{\ell^2(\Gamma)} +\mathcal E_{\mathrm{off}}$.
    \end{itemize}
\end{remark}

We call $H\colon D(H)\subset \ell^2(\Gamma)\to \ell^2(\Gamma)$, i.e., $H_{\mathrm{max}}$ on~\eqref{eq:VerySmallDomainFiniteInterRange}, the \emph{desired operator of finite interaction range} $m$ corresponding to $(H_{pq})_{p,q\in\Gamma}$.

Since $H$ has finite interaction range, only sites at
distance at most $m$ are coupled. As the diagonal part commutes with the localization
operators, only the off-diagonal contributions enter
the commutator. The off-diagonal edge energy
$\mathcal E_{\mathrm{off}}$ measures
these contributions individually, without cancellation.
The following theorem shows that controlling this
quantity by the graph norm of $H$ yields the required
averaged commutator estimate. By virtue of Remark~\ref{rem:PropertiesFiniteInteractionRangeUnbounded}, this is sufficient for Theorem~\ref{thm:window-pseudo-convergence} if $H=H_{\mathrm{max}}$.

We then specialize to weighted shifts and their finite
sums on $\Gamma=G$ (including $\mathbb Z^d$) and, more
generally, discrete groups with the above doubling property.

\begin{theorem}[Commutator property for finite-interaction-range operators]
\label{thm:finite-range-edge-energy}
Let $X=\Gamma \subset G$, and $H$ be the desired operator of finite interaction range $m\geq 0$ corresponding to a matrix $(H_{pq})_{p,q\in\Gamma}$ in $\C$.

Assume that there exists $C_E>0$ such that
\[
    \forall \psi\in D(H) : \mathcal E_{\mathrm{off}}(\psi)
    \leq
    C_E
    \bigl(
        \|H\psi\|_{\ell^2(\Gamma)}
        +
        \|\psi\|_{\ell^2(\Gamma)}
    \bigr).
\]
Then, $W_{L,x}D(H)\subset D(H)$ for all $L>0$ and $x\in G$. Moreover, with the
constant
\[
C_H \coloneqq 2\sqrt{N_m-1}C_D m C_E \geq 0,
\]
where\footnote{By Remark~\ref{rem:PropertiesFiniteInteractionRangeUnbounded}, we have a concrete bound on $N_m$, but the structure of $H$ may yield a smaller value, e.g., in the diagonal case, one has $N_m=1$, and hence $C_H=0$.} $N_m\coloneqq\sup_{p\in\Gamma}\#(B^\leq_m(p)\cap\Gamma)$,
\[
    \left(
        \int_G
        \bigl\|[H,W_{L,x}]\psi\bigr\|_{\ell^2(\Gamma)}^2
        \,\mathrm d\mu(x)
    \right)^{1/2}
    \leq
    \frac{C_H}{L}
    \|w_{L,e}\|_{L^2(G)}
    \bigl(
        \|H\psi\|_{\ell^2(\Gamma)}
        +
        \|\psi\|_{\ell^2(\Gamma)}
    \bigr)
\]
for all $L\geq m$ and all $\psi\in D(H)$.
\end{theorem}

\begin{proof}
We first establish the invariance of the domain under the localization
operators. Let $\psi\in D(H)$. Since $|w_{L,x}|\leq1$,
\[
    \|vW_{L,x}\psi\|_{\ell^2(\Gamma)}
    \leq
    \|v\psi\|_{\ell^2(\Gamma)}
\]
and
\[
    \mathcal E_{\mathrm{off}}(W_{L,x}\psi)^2
    =
    \sum_{p\in\Gamma}
    \sum_{\substack{q\in\Gamma\\0<d(p,q)\leq m}}
    |H_{pq}|^2|w_{L,x}(q)\psi(q)|^2
    \leq
    \mathcal E_{\mathrm{off}}(\psi)^2.
\]
Thus, $W_{L,x}\psi\in D(H)$.

For $\psi\in D(H)$, the diagonal part cancels in the commutator, and hence
\[
    ([H,W_{L,x}]\psi)(p)
    =
    \sum_{\substack{q\in\Gamma\\0<d(p,q)\leq m}}
    H_{pq}
    \bigl(
        w_{L,x}(q)-w_{L,x}(p)
    \bigr)
    \psi(q).
\]
Therefore, by Cauchy--Schwarz,
\[
    |([H,W_{L,x}]\psi)(p)|^2
    \leq
    (N_m-1)
    \sum_{\substack{q\in\Gamma\\0<d(p,q)\leq m}}
    |H_{pq}|^2
    |w_{L,x}(q)-w_{L,x}(p)|^2
    |\psi(q)|^2
\]
for $p\in\Gamma$.
Summing over $p$, integrating over $x$, and applying Tonelli's theorem gives
\[
\begin{aligned}
    \int_G
    \|[H,W_{L,x}]\psi\|_{\ell^2(\Gamma)}^2
    \,\mathrm d\mu(x)
    \leq{}&
    (N_m-1)
    \sum_{p\in\Gamma}
    \sum_{\substack{q\in\Gamma\\0<d(p,q)\leq m}}
    |H_{pq}|^2|\psi(q)|^2 \\
    &\qquad\qquad\times
    \int_G
    |w_{L,x}(q)-w_{L,x}(p)|^2
    \,\mathrm d\mu(x).
\end{aligned}
\]
Since the tent functions are $L^{-1}$-Lipschitz,
\[
    |w_{L,x}(q)-w_{L,x}(p)|
    \leq
    \frac{d(p,q)}{L}
    \leq
    \frac{m}{L}.
\]
Moreover, the integrand vanishes outside
$B^\leq_L(p)\cup B^\leq_L(q)\subset B^\leq_{L+m}(q)$. For $L\geq m$, the doubling property and left-invariance give
\[
    \mu(B^\leq_{L+m}(q))
    \leq \mu(B^\leq_{2L}(e))
    \leq C_D^2\mu(B^\leq_{L/2}(e)).
\]
Since $w_{L,e}\geq \frac12$ on $B^\leq_{L/2}(e)$,
\[
    \mu(B^\leq_{L/2}(e))
    \leq 4\|w_{L,e}\|_{L^2(G)}^2.
\]
Consequently,
\[
    \int_G |w_{L,x}(q)-w_{L,x}(p)|^2\,\mathrm d\mu(x)
    \leq
    \frac{4C_D^2m^2}{L^2}
    \|w_{L,e}\|_{L^2(G)}^2.
\]
Substituting this estimate into the previous bound yields
\[
    \int_G
    \|[H,W_{L,x}]\psi\|_{\ell^2(\Gamma)}^2
    \,\mathrm d\mu(x)
    \leq
    \frac{4(N_m-1)C_D^2m^2}{L^2}
    \|w_{L,e}\|_{L^2(G)}^2
    \mathcal E_{\mathrm{off}}(\psi)^2.
\]
The assumed graph-norm bound for $\mathcal E_{\mathrm{off}}$ now yields
\[
    \int_G
    \|[H,W_{L,x}]\psi\|_{\ell^2(\Gamma)}^2
    \,\mathrm d\mu(x)
    \leq
    \frac{4(N_m-1)C_D^2m^2C_E^2}{L^2}
    \|w_{L,e}\|_{L^2(G)}^2
    \bigl(
        \|H\psi\|_{\ell^2(\Gamma)}
        +
        \|\psi\|_{\ell^2(\Gamma)}
    \bigr)^2.
\]
Taking square roots proves the claim.
\end{proof}

A single weighted shift is easy to treat.
Its off-diagonal edge energy coincides with $\|H\psi\|_2$
whenever the shift is nontrivial. Thus, no additional
restriction on the growth of the weight is needed,
and the commutator estimate involves only $\|H\psi\|_2$.

\begin{lemma}[Commutator property for a weighted shift]
\label{lem:weighted-shift-graph-comm}
Let $X=\Gamma = G$ be a discrete group and fix $s\in\Gamma$. Define the shift operator in direction $s$ as 
\[
    (S_s\psi)(p)\coloneqq \psi(ps),
    \qquad p\in\Gamma,
\]
and let $a\colon\Gamma\to\C$ be arbitrary, possibly unbounded. For the corresponding multiplication operator $D_a$, consider
$D_aS_s$ on its domain
\[
    D(D_aS_s)\coloneqq
    \left\{
        \psi\in\ell^2(\Gamma):
        p\mapsto a(p)\psi(ps)\in\ell^2(\Gamma)
    \right\}.
\]
Then, $D_aS_s=H$, where $H$ is the desired operator of finite interaction range $d(e,s)$ corresponding to the matrix $(a(p)\delta_{ps}(q))_{p,q\in\Gamma}$.
Moreover, setting $C_s\coloneqq 2C_Dd(e,s)$, we have,
for all sufficiently large $L$ and all $\psi\in D(H)$,
\[
    \left(
        \int_\Gamma
        \|[H,W_{L,x}]\psi\|_{\ell^2(\Gamma)}^2
        \,\mathrm d\mu(x)
    \right)^{1/2}
    \leq
    \frac{C_s}{L}
    \|w_{L,e}\|_{L^2(\Gamma)}
    \|H\psi\|_{\ell^2(\Gamma)}.
\]
\end{lemma}

\begin{proof}
$D_aS_s=H$ is clear.
If $s=e$, then $H=D_a$ is diagonal and $[H,W_{L,x}]=0$, so the bound is immediate. Let $s\neq e$. Then $H$ has
interaction range $m\coloneqq d(e,s)$ and vanishing diagonal with
\[
    \mathcal E_{\mathrm{off}}(\psi)^2
    =
    \sum_{p\in\Gamma}|a(p)|^2|\psi(ps)|^2
    =
    \|H\psi\|_{\ell^2(\Gamma)}^2.
\]
Thus, we can imitate Theorem~\ref{thm:finite-range-edge-energy} with
$C_E=1$ and without the necessity of an $N_m$, and its edge-energy estimate yields
the assertion.
\end{proof}

For a finite sum of weighted shifts, the individual
contributions may cancel, so their norms need not be
controlled by the norm of their sum. The following
assumption supplies this control in the graph norm of $H$.
It is precisely the edge-energy condition from
Theorem~\ref{thm:finite-range-edge-energy}, expressed
in terms of the individual shifts.

\begin{theorem}[Commutator property for finite sums of weighted shifts]
\label{thm:unbounded-finite-range-group}
Let $X=\Gamma = G$ be a discrete group, and let $F\subset\Gamma$ be finite with $F_0\coloneqq F\setminus\{e\}\neq\varnothing$. For $s\in F$, let $a_s:\Gamma\to\C$ be arbitrary, possibly unbounded, and consider $D_{a_s}S_s$ as in Lemma~\ref{lem:weighted-shift-graph-comm}.
Consider
\[
    \sum_{s\in F}H_s
    \qquad\text{ with domain }
    \bigcap_{s\in F}D(H_s).
\]
Then, this is the desired operator $H$ of finite interaction range $m_F\coloneqq\max_{s\in F_0}d(e,s)$ corresponding to the matrix $(\sum_{s\in F}a_s(p)\delta_{ps}(q))_{p,q\in\Gamma}$.
 Assume that there exists $C_{\mathrm{br}}>0$ such that
\[
    \left(
        \sum_{s\in F_0}
        \|H_s\psi\|_{\ell^2(\Gamma)}^2
    \right)^{1/2}
    \leq
    C_{\mathrm{br}}
    \bigl(
        \|H\psi\|_{\ell^2(\Gamma)}
        +
        \|\psi\|_{\ell^2(\Gamma)}
    \bigr),
    \qquad \psi\in D(H).
\]
Then $W_{L,x}D(H)\subset D(H)$ for all $L>0$ and $x\in\Gamma$. Moreover, for all sufficiently large $L$ and all $\psi\in D(H)$,
\[
    \left(
        \int_\Gamma
        \|[H,W_{L,x}]\psi\|_{\ell^2(\Gamma)}^2
        \,\mathrm d\mu(x)
    \right)^{1/2}
    \leq
    \frac{2\sqrt{|F_0|}\,C_Dm_FC_{\mathrm{br}}}{L}
    \|w_{L,e}\|_{L^2(\Gamma)}
    \bigl(
        \|H\psi\|_{\ell^2(\Gamma)}
        +
        \|\psi\|_{\ell^2(\Gamma)}
    \bigr).
\]
\end{theorem}

\begin{proof}

Since the points $ps$,
$s\in F_0$, are distinct,
\[
    \mathcal E_{\mathrm{off}}(\psi)^2
    =
    \sum_{s\in F_0}\sum_{p\in\Gamma}
    |a_s(p)|^2|\psi(ps)|^2
    =
    \sum_{s\in F_0}\|H_s\psi\|_{\ell^2(\Gamma)}^2,
\]
and $H=\sum_{s\in F}H_s$ is immediate.
Thus, imitating Theorem~\ref{thm:finite-range-edge-energy}
with
$C_E=C_{\mathrm{br}}$ and
$|F_0|$ as a replacement for $N_{m_F}-1$ yields the assertion.
\end{proof}

To make the edge-energy condition more concrete, we now consider band operators on $\ell^2(\mathbb Z)$.
In this setting, the weighted shifts describe the
off-diagonals. Boundedness or suitable dominance relations between their coefficients provide explicit
sufficient conditions for the required graph-norm bound and for $H=H_{\mathrm{max}}$.

\begin{example}[Unbounded band operators on $\ell^2(\mathbb Z)$]
Let $X=\Gamma=G=\mathbb Z$ with its usual metric, and consider
\[
    (H\psi)(k)
    \coloneqq
    v(k)\psi(k)
    +
    \sum_{j\in J}a_j(k)\psi(k+j),
\]
where $k\in\mathbb{Z}$, $J\subset\mathbb Z\setminus\{0\}$ is finite and non-empty,
and $v,a_j\colon\mathbb Z\to\mathbb C$ may be unbounded. Write
$H_j\coloneqq D_{a_j}S_j$ and equip $H$ with the domain
\[
    D(H)\coloneqq
    \left\{
        \psi\in\ell^2(\mathbb Z):
        v\psi\in\ell^2(\mathbb Z),\;
        H_j\psi\in\ell^2(\mathbb Z)\ \text{for all }j\in J
    \right\}.
\]
Here $v$ describes the main diagonal, while $a_j$ describes
the diagonal with offset $j$. It is easy to see that these are actually desired operators of finite interaction range in the regime of Theorem~\ref{thm:unbounded-finite-range-group}.
Set $n\coloneqq |J|$ and
\[
    b(k)\coloneqq
    \Bigl(\sum_{j\in J}|a_j(k-j)|^2\Bigr)^{1/2}.
\]
Changing variables gives
\[
    \mathcal E_{\mathrm{off}}(\psi)^2
    =
    \sum_{j\in J}\|H_j\psi\|_2^2
    =
    \|b\psi\|_2^2.
\]
Each of the following conditions is sufficient for the
branch-energy assumption of
Theorem~\ref{thm:unbounded-finite-range-group} and implies $H=H_{\mathrm{max}}$ by virtue of~\eqref{eq:HequalsHmaxequalsHmin}. Therefore, we can apply Theorem~\ref{thm:window-pseudo-convergence}.

In the following, we always assume $\psi\in \ell^2(\mathbb Z)$.
\begin{enumerate}
\item \emph{Dominance of the main diagonal.}
There exist $\beta\geq0$ and $0\leq\eta<n^{-1/2}$ such that
\[
    b(k)\leq\eta|v(k)|+\beta,
    \qquad k\in\mathbb Z.
\]
Since
\[
    \mathcal E_{\mathrm{off}}(\psi)
    =\|b\psi\|_2
    \leq\eta\|v\psi\|_2+\beta\|\psi\|_2
    \leq\eta\|H\psi\|_2
    +\eta\sqrt n\,\mathcal E_{\mathrm{off}}(\psi)
    +\beta\|\psi\|_2,
\]
absorption gives
\[
    \mathcal E_{\mathrm{off}}(\psi)
    \leq
    \frac{\eta}{1-\eta\sqrt n}\|H\psi\|_2
    +
    \frac{\beta}{1-\eta\sqrt n}\|\psi\|_2.
\]

\item \emph{Individual comparison with the main diagonal.}
There exist $\eta_j,\beta_j\geq0$ such that
\[
    |a_j(k-j)|\leq\eta_j|v(k)|+\beta_j,
    \qquad k\in\mathbb Z,\quad j\in J,
\]
where $\eta\coloneqq\sum_{j\in J}\eta_j<1$.
Setting $\beta:=\sum_{j\in J}\beta_j$ and
$\mathcal A_{\mathrm{off}}(\psi)\coloneqq\sum_{j\in J}\|H_j\psi\|_2$,
we obtain
\[
    \mathcal A_{\mathrm{off}}(\psi)
    \leq\eta\|v\psi\|_2+\beta\|\psi\|_2
    \leq\eta\|H\psi\|_2
    +\eta\mathcal A_{\mathrm{off}}(\psi)
    +\beta\|\psi\|_2.
\]
Thus
\[
    \mathcal E_{\mathrm{off}}(\psi)
    \leq\mathcal A_{\mathrm{off}}(\psi)
    \leq
    \frac{\eta}{1-\eta}\|H\psi\|_2
    +
    \frac{\beta}{1-\eta}\|\psi\|_2.
\]

\item \emph{Bounded off-diagonals.}
All $a_j$ are bounded. Then $b$ is bounded and
\[
    \mathcal E_{\mathrm{off}}(\psi)
    \leq\|b\|_\infty\|\psi\|_2,
\]
without any growth restriction on $v$.

\item \emph{Dominance of one off-diagonal.}
The main diagonal $v$ is bounded, and there exist $j_0\in J$
and $\eta_j,\beta_j\geq0$ such that
\[
    |a_j(k-j)|
    \leq\eta_j|a_{j_0}(k-j_0)|+\beta_j,
    \qquad k\in\mathbb Z,\quad j\in J\setminus\{j_0\},
\]
where
\[
    \eta\coloneqq\sum_{j\in J\setminus\{j_0\}}\eta_j<1,
    \qquad
    \beta\coloneqq\sum_{j\in J\setminus\{j_0\}}\beta_j.
\]
Then
\[
    \sum_{j\neq j_0}\|H_j\psi\|_2
    \leq\eta\|H_{j_0}\psi\|_2+\beta\|\psi\|_2.
\]
The reverse triangle inequality therefore gives
\[
    (1-\eta)\|H_{j_0}\psi\|_2
    \leq\|H\psi\|_2+(\|v\|_\infty+\beta)\|\psi\|_2,
\]
and hence
\[
    \mathcal E_{\mathrm{off}}(\psi)
    \leq(1+\eta)\|H_{j_0}\psi\|_2+\beta\|\psi\|_2
    \leq
    \frac{1+\eta}{1-\eta}\|H\psi\|_2
    +
    \frac{(1+\eta)\|v\|_\infty+2\beta}{1-\eta}\|\psi\|_2.
\]
\end{enumerate}
\end{example}

%% file: schroedinger_example.tex
\section{Schr\"odinger operators}
\label{sec:schrodinger-core-closure}

In this section we show that Theorem \ref{thm:window-pseudo-convergence} on pseudospectral bounds and their convergence applies to Schr\"odinger operators $H$ on $L^2(\mathbb R^d)$ with possibly unbounded potentials, when we use a $C^\infty_c$ family of localization functions $(\chi_{L,x})_{L>0,\,x\in G}$. 

So, let $X=G=\R^d$ with $d\in\N$. In Subsections 6.1 and 6.2 we show that, and under which concrete assumptions on the potential and the localization profile $\eta$, $H$ and the family $(\chi_{L,x})_{L>0,\,x\in G}$ have the commutator property from Definition \ref{def:commute}.

In Subsection 6.3 we then argue why this already yields the same results for $H^*$, i.e., why Theorem~\ref{thm:window-pseudo-convergence} is applicable.

We first carry out the calculations on compactly supported
smooth functions and then extend the result to the closure of the operator.

Let $V\in L^\infty_{\mathrm{loc}}(\mathbb R^d)$ be complex-valued and
satisfy $\operatorname{Re}V\geq -V_0$ for some $V_0\geq0$. Define
\begin{equation}\label{eq:definSchrodingerOnCcinf}
    D(H_0)\coloneqq C_c^\infty(\mathbb R^d),
    \qquad H_0\psi\coloneqq -\Delta\psi+V\psi.
\end{equation}
The operator $H_0$ is densely defined and closable: its adjoint has $C_c^\infty(\mathbb R^d)$, where it acts as
$-\Delta+V^*$, in its domain. We set $H\coloneqq\overline{H_0}$.

All norms below are $L^2$-norms unless indicated otherwise, and the
inner product is linear in its first argument. Fix a nonzero real-valued localization profile $\eta\in C_c^\infty(B^\leq_1(0))$ and set
\[
    \chi_{L,x}(y)\coloneqq\eta\bigl((y-x)/L\bigr),
    \qquad L\geq1,\quad x,y\in\mathbb R^d.
\]
We use the same notation for these functions and the corresponding
multiplication operators.

\subsection{Estimates on compactly supported smooth functions}

\begin{lemma}\label{lem:schrodinger-core-estimates}
For every $\psi\in C_c^\infty(\mathbb R^d)$,
\begin{equation}\label{eq:schrodinger-core-gradient}
    \|\nabla\psi\|_2
    \leq \frac12\|H_0\psi\|_2
    +\left(\frac12+\sqrt{V_0}\right)\|\psi\|_2.
\end{equation}
Moreover, for every $L\geq1$ and $x\in\mathbb R^d$, multiplication by $\chi_{L,x}$ preserves
$C_c^\infty(\mathbb R^d)$, and
\begin{equation}\label{eq:schrodinger-core-average}
    \left(\int_{\mathbb R^d}
    \|[H_0,\chi_{L,x}]\psi\|_2^2\,\mathrm dx\right)^{1/2}
    \leq \frac{C_H}{L}\|\chi_{L,0}\|_2
    \bigl(\|H_0\psi\|_2+\|\psi\|_2\bigr),
\end{equation}
where
\begin{equation*}%
    C_H\coloneqq\left(1+2\sqrt{V_0}\right)
    \frac{\|\nabla\eta\|_2}{\|\eta\|_2}
    +\frac{\|\Delta\eta\|_2}{\|\eta\|_2}.
\end{equation*}
\end{lemma}

\begin{proof}
Integration by parts and the lower bound on $\operatorname{Re}V$ give
\[
    \|\nabla\psi\|_2^2
    \leq \operatorname{Re}\langle H_0\psi,\psi\rangle+V_0\|\psi\|_2^2
    \leq \|H_0\psi\|_2\|\psi\|_2+V_0\|\psi\|_2^2.
\]
Taking square roots and using
$\sqrt{a+b}\leq\sqrt a+\sqrt b$, followed by
$\sqrt{cd}\leq(c+d)/2$, gives
\eqref{eq:schrodinger-core-gradient}.

Since $\chi_{L,x}\psi\in C_c^\infty(\mathbb R^d)$ and
multiplication by $V$ commutes with multiplication by $\chi_{L,x}$,
the product rule yields
\begin{equation}\label{eq:schrodinger-core-commutator}
    [H_0,\chi_{L,x}]\psi
    =-2\nabla\chi_{L,x}\cdot\nabla\psi
    -(\Delta\chi_{L,x})\psi.
\end{equation}
For fixed $y\in\mathbb R^d$, a change of variables gives
\begin{align}
    \int_{\mathbb R^d}|\nabla\chi_{L,x}(y)|^2\,\mathrm dx
    &=L^{d-2}\|\nabla\eta\|_2^2,
    \label{eq:schrodinger-gradient-scaling}\\
    \int_{\mathbb R^d}|\Delta\chi_{L,x}(y)|^2\,\mathrm dx
    &=L^{d-4}\|\Delta\eta\|_2^2.
    \label{eq:schrodinger-laplacian-scaling}
\end{align}

We now insert \eqref{eq:schrodinger-core-commutator} into the
averaged norm. The triangle inequality in
$L^2(\mathbb R^d_x\times\mathbb R^d_y)$ gives
\[
\begin{aligned}
\left(\int_{\mathbb R^d}
    \|[H_0,\chi_{L,x}]\psi\|_2^2\,\mathrm dx\right)^{1/2} 
&\quad=
\left(\int_{\mathbb R^d}\int_{\mathbb R^d}
\left|2\nabla\chi_{L,x}(y)\cdot\nabla\psi(y)
+(\Delta\chi_{L,x}(y))\psi(y)\right|^2
\,\mathrm dy\,\mathrm dx\right)^{1/2}\\
&\quad\leq
2\left(\int_{\mathbb R^d}\int_{\mathbb R^d}
|\nabla\chi_{L,x}(y)|^2|\nabla\psi(y)|^2
\,\mathrm dy\,\mathrm dx\right)^{1/2}\\
&\quad +
\left(\int_{\mathbb R^d}\int_{\mathbb R^d}
|\Delta\chi_{L,x}(y)|^2|\psi(y)|^2
\,\mathrm dy\,\mathrm dx\right)^{1/2}.
\end{aligned}
\]
By Tonelli's theorem, we can integrate first with respect to $x$.
Using the two scaling identities \eqref{eq:schrodinger-gradient-scaling} and \eqref{eq:schrodinger-laplacian-scaling}, we obtain
\[
\begin{aligned}
\left(\int_{\mathbb R^d}
    \|[H_0,\chi_{L,x}]\psi\|_2^2\,\mathrm dx\right)^{1/2}
&\quad\leq
2\left(\int_{\mathbb R^d}|\nabla\psi(y)|^2
    \int_{\mathbb R^d}|\nabla\chi_{L,x}(y)|^2
    \,\mathrm dx\,\mathrm dy\right)^{1/2}\\
&\quad+
\left(\int_{\mathbb R^d}|\psi(y)|^2
    \int_{\mathbb R^d}|\Delta\chi_{L,x}(y)|^2
    \,\mathrm dx\,\mathrm dy\right)^{1/2}\\
&\quad=
2L^{d/2-1}\|\nabla\eta\|_2\|\nabla\psi\|_2
+L^{d/2-2}\|\Delta\eta\|_2\|\psi\|_2.
\end{aligned}
\]
Combining the preceding estimate with
\eqref{eq:schrodinger-core-gradient} gives
\[
\begin{aligned}
\left(\int_{\mathbb R^d}
    \|[H_0,\chi_{L,x}]\psi\|_2^2\,\mathrm dx\right)^{1/2}
&\leq L^{d/2-1}\|\nabla\eta\|_2\|H_0\psi\|_2\\
&\quad+\left((1+2\sqrt{V_0})L^{d/2-1}\|\nabla\eta\|_2
          +L^{d/2-2}\|\Delta\eta\|_2\right)\|\psi\|_2.
\end{aligned}
\]
Since $\|\chi_{L,0}\|_2=L^{d/2}\|\eta\|_2$, this can be written as
\[
\begin{aligned}
\left(\int_{\mathbb R^d}
    \|[H_0,\chi_{L,x}]\psi\|_2^2\,\mathrm dx\right)^{1/2}
&\leq \frac{\|\chi_{L,0}\|_2}{L}
\biggl[
\frac{\|\nabla\eta\|_2}{\|\eta\|_2}\|H_0\psi\|_2\\
&\qquad+
\left((1+2\sqrt{V_0})\frac{\|\nabla\eta\|_2}{\|\eta\|_2}
+\frac{1}{L}\frac{\|\Delta\eta\|_2}{\|\eta\|_2}\right)
\|\psi\|_2
\biggr].
\end{aligned}
\]
As $L\geq1$, we may choose
\[
C_H=(1+2\sqrt{V_0})\frac{\|\nabla\eta\|_2}{\|\eta\|_2}
    +\frac{\|\Delta\eta\|_2}{\|\eta\|_2}.
\]
Therefore,
\[
\left(\int_{\mathbb R^d}
    \|[H_0,\chi_{L,x}]\psi\|_2^2\,\mathrm dx\right)^{1/2}
\leq \frac{C_H}{L}\|\chi_{L,0}\|_2
    \bigl(\|H_0\psi\|_2+\|\psi\|_2\bigr).
\]
\end{proof}

\subsection{Commutator estimate for the closed Schr\"odinger operator}
\mbox{}\\
By the definition of $H=\overline{H_0}$, for every
$\psi\in D(H)$, there exist
$\psi_n\in C_c^\infty(\mathbb R^d)$ such that
\[
    \psi_n\longrightarrow\psi,
    \qquad H_0\psi_n\longrightarrow H\psi
    \quad\text{in }L^2(\mathbb R^d).
\]
Applying \eqref{eq:schrodinger-core-gradient} to
$\psi_n-\psi_m$ shows that $(\psi_n)$ is Cauchy in
$H^1(\mathbb R^d)$. Thus $\psi_n\to\psi$ in $H^1$.
We use this approximation to extend the cutoff invariance and the
commutator estimate from $H_0$ to $H$.

\begin{theorem}\label{thm:schrodinger-closed-commutator}
For every $L\geq1$ and $x\in\mathbb R^d$, multiplication by $\chi_{L,x}\in C_c^\infty(\mathbb R^d)$ preserves $D(H)$,
and
\begin{equation}\label{eq:schrodinger-closed-commutator}
    [H,\chi_{L,x}]\psi=-2\nabla\chi_{L,x}\cdot\nabla\psi-(\Delta\chi_{L,x})\psi,
    \qquad \psi\in D(H).
\end{equation}
For every $L\geq1$ and $\psi\in D(H)$,
\begin{equation*}%
    \left(\int_{\mathbb R^d}
    \|[H,\chi_{L,x}]\psi\|_2^2\,\mathrm dx\right)^{1/2}
    \leq \frac{C_H}{L}\|\chi_{L,0}\|_2
    \bigl(\|H\psi\|_2+\|\psi\|_2\bigr),
\end{equation*}
with the same constant $C_H$ as in
Lemma~\ref{lem:schrodinger-core-estimates}.
\end{theorem}

\begin{proof}
Fix $\psi\in D(H)$ and choose $\psi_n\in C_c^\infty(\mathbb R^d)$, so that
$\psi_n\to\psi$ in $H^1$ and $H_0\psi_n\to H\psi$ in $L^2$.
For fixed $L\geq1$ and $x\in\mathbb R^d$, we have
$\chi_{L,x}\psi_n\in C_c^\infty(\mathbb R^d)$, and
\[
H_0(\chi_{L,x}\psi_n)
=\chi_{L,x}H_0\psi_n
 -2\nabla\chi_{L,x}\cdot\nabla\psi_n
 -(\Delta\chi_{L,x})\psi_n.
\]
All terms on the right converge in $L^2$, while
$\chi_{L,x}\psi_n\to\chi_{L,x}\psi$ in $L^2$. Since $H$ is closed,
it follows that $\chi_{L,x}\psi\in D(H)$ and
\[
H(\chi_{L,x}\psi)
=\chi_{L,x}H\psi
 -2\nabla\chi_{L,x}\cdot\nabla\psi
 -(\Delta\chi_{L,x})\psi.
\]
This proves \eqref{eq:schrodinger-closed-commutator}. The same
convergence, together with $H_0\psi_n\to H\psi$, gives
\[
[H_0,\chi_{L,x}]\psi_n
\longrightarrow [H,\chi_{L,x}]\psi
\quad\text{in }L^2
\]
for every fixed $x$. Fatou's lemma and
\eqref{eq:schrodinger-core-average} now yield
\[
\begin{aligned}
\left(\int_{\mathbb R^d}
    \|[H,\chi_{L,x}]\psi\|_2^2\,\mathrm dx\right)^{1/2}
&\leq \liminf_{n\to\infty}
    \left(\int_{\mathbb R^d}
    \|[H_0,\chi_{L,x}]\psi_n\|_2^2\,\mathrm dx\right)^{1/2}\\
&\leq \frac{C_H}{L}\|\chi_{L,0}\|_2
    \bigl(\|H\psi\|_2+\|\psi\|_2\bigr).
\end{aligned}
\]
\end{proof}

\subsection{Equality of minimal and maximal domain, $H_{min} = H_{max}$}

\begin{proposition}\label{prop:schrodinger-min-max}
Let
\[
    V\in L^2_{\mathrm{loc}}(\mathbb R^d;\mathbb C),
    \qquad \operatorname{Re}V\geq -V_0,
    \qquad V_0\geq0,
\]
and define
\[
    H_0u\coloneqq -\Delta u+Vu,
    \qquad D(H_0)\coloneqq C_c^\infty(\mathbb R^d).
\]
Then $H_0$ is closable, and its closure $H_{\min}\coloneqq\overline{H_0}$
coincides with the maximal (distributional) operator
\[
    H_{\max}u\coloneqq -\Delta u+Vu,
    \qquad
    D(H_{\max})
    \coloneqq
    \left\{
        u\in L^2(\mathbb R^d):
        -\Delta u+Vu\in L^2(\mathbb R^d)
        \text{ distributionally}
    \right\}.
\]
\end{proposition}

\begin{proof}
We adapt the proof strategy of
\cite[Theorem~X.28]{ReedSimonII} to complex potentials.

\smallskip
\noindent\emph{1. Minimal and maximal operators.}
For every $u\in L^2$, Hölder's inequality gives
$Vu\in L^1_{\mathrm{loc}}$, so the distributional expression
$-\Delta u+Vu$ is well defined. For $W\in\{V,V^*\}$, we write $H_0(W)$ for $-\Delta+W$
on $C_c^\infty(\mathbb R^d)$, $H_{\min}(W)$ for its closure,
and $H_{\max}(W)$ for the corresponding maximal distributional
realization.

The definition of the adjoint then gives
\[
    H_0(V)^*=H_{\max}(V^*),
    \qquad
    H_0(V^*)^*=H_{\max}(V).
\]
Since $C_c^\infty\subset D(H_{\max}(V^*))$, the adjoint
$H_0(V)^*$ is densely defined, and hence $H_0$ is closable.
Moreover, $H_{\max}$ is closed and extends $H_0$, so
\[
    H_{\min}\subset H_{\max},
    \qquad
    H_{\min}^*=H_{\max}(V^*).
\]

\smallskip
\noindent\emph{2. Reduction to injectivity.}
Fix $a>V_0$ and put $c:=a-V_0>0$. We first show that
$H_{\min}=H_{\max}$ follows from
\begin{equation}\label{eq:kernel-condition}
    \ker\bigl(H_{\max}(V^*)+a\bigr)=\{0\}.
\end{equation}
The standard range--kernel identity gives
\begin{equation}\label{eq:range-kernel}
    \operatorname{Ran}(H_{\min}+a)^\perp
    =
    \ker\bigl((H_{\min}+a)^*\bigr)
    =
    \ker\bigl(H_{\max}(V^*)+a\bigr).
\end{equation}
By \eqref{eq:range-kernel} and \eqref{eq:kernel-condition},
the range of $H_{\min}+a$ is dense in $L^2(\mathbb R^d)$.

For $\varphi\in C_c^\infty$, integration by parts gives
\[
    \operatorname{Re}
    \langle (H_0+a)\varphi,\varphi\rangle
    =
    \|\nabla\varphi\|_2^2
    +
    \int_{\mathbb R^d}
        (\operatorname{Re}V+a)|\varphi|^2
    \geq
    c\|\varphi\|_2^2.
\]
By Cauchy--Schwarz and passage to the closure, it follows that
\[
    \|(H_{\min}+a)u\|_2
    \geq c\|u\|_2,
    \qquad u\in D(H_{\min}).
\]
Since $H_{\min}+a$ is closed, this estimate implies that its range
is closed. Consequently,
\[
    \operatorname{Ran}(H_{\min}+a)=L^2(\mathbb R^d).
\]

Since $a$ is real, complex conjugation maps
$\ker(H_{\max}(V)+a)$ onto
$\ker(H_{\max}(V^*)+a)$. Thus
\eqref{eq:kernel-condition} also gives
\begin{equation}\label{eq:kernel-condition-V}
    \ker(H_{\max}+a)=\{0\}.
\end{equation}
Now let $v\in D(H_{\max})$. By surjectivity, there exists
$u\in D(H_{\min})$ such that
\begin{equation*}%
    (H_{\min}+a)u=(H_{\max}+a)v.
\end{equation*}
With $H_{\min}\subset H_{\max}$, we obtain
\begin{equation*}%
    (H_{\max}+a)(v-u)=0.
\end{equation*}
By the triviality of the kernel in \eqref{eq:kernel-condition-V},
this implies
$v=u\in D(H_{\min})$, proving the reverse inclusion.

\smallskip
\noindent\emph{3. Kato's inequality and mollification.}
Let
\[
    w\in\ker\bigl(H_{\max}(V^*)+a\bigr).
\]
Then
\[
    \Delta w=(V^*+a)w
\]
in the sense of distributions. Since
$V^*+a,w\in L^2_{\mathrm{loc}}$, H\"older's inequality
gives $(V^*+a)w\in L^1_{\mathrm{loc}}$. Therefore
$w,\Delta w\in L^1_{\mathrm{loc}}$, and the distributional Kato
inequality \cite[Theorem~X.27]{ReedSimonII} applies:
\[
    \Delta|w|
    \geq
    \operatorname{Re}\left(
        \frac{\overline w}{|w|}\,\Delta w
    \right)
    =
    (\operatorname{Re}V+a)|w|
    \geq c|w|.
\]
Here the quotient is defined to be zero on $\{w=0\}$, and the
inequalities are understood distributionally.

Following the mollification argument in the proof of
\cite[Theorem~X.28]{ReedSimonII}, choose
$\rho\in C_c^\infty(\mathbb R^d)$ with $\rho\geq0$ and
$\int_{\mathbb R^d}\rho=1$, and set
\[
    \rho_\delta(x):=\delta^{-d}\rho(x/\delta),
    \qquad
    z_\delta:=|w|*\rho_\delta,
    \qquad \delta>0.
\]
Then $z_\delta\geq0$, and Young's inequality gives
\[
    \partial^\alpha z_\delta
    =
    |w|*(\partial^\alpha\rho_\delta)\in L^2,
    \qquad |\alpha|\leq2.
\]
In particular, $z_\delta\in C^\infty\cap H^2$.

Convolving the distributional inequality $\Delta|w|\geq0$
with the nonnegative mollifier $\rho_\delta$ gives
\begin{equation*}%
    \Delta z_\delta
    =
    (\Delta|w|)*\rho_\delta
    \geq0
\end{equation*}
pointwise on $\mathbb R^d$.
Since $z_\delta\geq0$ and $z_\delta\in H^2$, integration by parts
yields
\begin{equation*}%
    0
    \leq
    \int_{\mathbb R^d}z_\delta\Delta z_\delta
    =
    -\|\nabla z_\delta\|_2^2
    \leq0.
\end{equation*}
Hence $z_\delta$ is constant and therefore zero, since
$z_\delta\in L^2(\mathbb R^d)$. Finally, $z_\delta\to|w|$ in $L^2$ as
$\delta\to0$, and therefore $w=0$.
\end{proof}
Considering our original Schr\"odinger operator $H$, i.e., the closure of~\eqref{eq:definSchrodingerOnCcinf}, Proposition~\ref{prop:schrodinger-min-max} proves that the adjoint $H^*$ is once again the closure of~\eqref{eq:definSchrodingerOnCcinf}, where we replaced $V$ with $V^*$. Hence, Theorem~\ref{thm:schrodinger-closed-commutator} holds for both $H$ and $H^*$ with the same constant $C_H$, and we can apply Theorem~\ref{thm:window-pseudo-convergence}.